\documentclass[12pt,reqno]{amsart}
\usepackage{amsmath,amsthm,amscd, amssymb}

\swapnumbers

\newcommand{\qlbarstar}{\bar{\mathbf Q}_l^{\times}}

\renewcommand{\l}{\lambda}
\renewcommand{\L}{ \Lambda}
 \newcommand{\s}{{\sigma}}
\renewcommand{\det}{\operatorname{det}}

\newcommand{\Hom}{\operatorname{Hom}}
\newcommand{\Lie}{\operatorname{Lie}}

\newcommand{\Qbar}{\bar {\mathbf Q}}

\newcommand{\Gal}{\operatorname{Gal}}

\newcommand{\G}{\Gamma}

\newcommand{\vl}{\operatorname{vol}}

\newcommand{\Fbar}{\overline F}

\newcommand{\Ad}{\operatorname{Ad}}

\newcommand{\cok}{\operatorname{cok}}

\newcommand{\Ext}{\operatorname{Ext}}

\newcommand{\Lbar}{\overline L}

\newtheorem{thm}[subsection]{Theorem}
\newtheorem{cor}[subsection]{Corollary}
\newtheorem{lem}[subsection]{Lemma}
\newtheorem{prop}[subsection]{Proposition}

\theoremstyle{definition}

\newcommand{\bS}{\mathbf S}
\newcommand{\bC}{\mathbb C}
\newcommand{\val}{\operatorname{val}}
\newcommand{\SL}{\text{\rm SL}(2)}
\newcommand{\quash}[1]{}    
\renewcommand{\mod}[1]{\ (\mathop{\rm mod} #1)}
\renewcommand{\epsilon}{\varepsilon}
\newcommand{\parens}[1]{{\rm{(}}#1{\rm{)}}} 
\newcommand{\abs}[1]{\left| #1 \right|}

\numberwithin{equation}{subsection}

\begin{document}

\title{Homology of affine Springer fibers in the unramified case}
\author{Mark Goresky${}^1$}\thanks{1.  partially supported by N.S.F. grant no.
DMS-9900324 and DMS-0139986}
\author {Robert Kottwitz${}^2$}\thanks{2.  partially supported by N.S.F. grant
no. DMS-0071971}
\author {Robert MacPherson}
\keywords{loop Grassmannian, affine Springer fiber, fundamental lemma,
equivariant homology}
\subjclass{primary: 22E67,55N91 secondary: 14M15,11D88, 11F85}
\begin{abstract}
   Assuming a certain ``purity" conjecture, we derive a formula
for the (complex) cohomology groups of the affine Springer fiber
corresponding to any unramified regular semi-simple element.
We use this calculation to present a complex analog of the
fundamental lemma for function fields.  We show that the ``kappa"
orbital integral which arises in the fundamental lemma is equal
to the Lefschetz trace of Frobenius acting on the \'etale
cohomology of a related variety.
\end{abstract}
\maketitle

\section{Introduction}  Let $F = \mathbb C((\epsilon))$ be the field of formal
Laurent series and let $\mathfrak o = \mathbb C[[\epsilon]]$ be its integer
ring of formal power series. A fundamental object of study in modern
representation theory is the affine Grassmannian $X = G\left(F \right)/ G
\left(  \mathfrak o \right)$ associated to a complex reductive algebraic group
$G.$  Let $\mathfrak g = \text{Lie}(G).$  Each regular semisimple element
$\gamma \in \mathfrak g(\mathfrak o)$ determines a ``vector field' on $X$
whose fixed point set
\[ X_{\gamma} = \left\{ xG(\mathfrak o) \in X:\ \text{Ad}(x^{-1})(\gamma) \in
\mathfrak g(\mathfrak o) \right\} \]
is a finite dimensional increasing union of complex projective varieties,
which is known as the affine Springer fiber corresponding to $\gamma.$  The
study of affine Springer fibers was initiated in \cite{KL}, where many of the
basic properties of these varieties were described.

In this paper we suppose $T \subset G$ is a maximal torus defined over
$\mathbb C$ and we consider the affine Springer fiber $X_{\gamma}$ associated
to regular semisimple elements $\gamma \in \mathfrak t(\mathfrak o),$ where
$\mathfrak t = \text{Lie}(T).$  (Not all affine Springer fibers are of this
type:  since $T$ is split over $F$ we refer to this as the unramified case.)
The complex torus $T(\mathbb C)$ acts on $X_{\gamma}.$  In Theorem
\ref{thm-cohomology}, the first principal result of this paper, we give an
explicit description of the $T(\mathbb C)-$equivariant homology of
$X_{\gamma}$ under the assumption that the ordinary homology of $X_{\gamma}$
is pure, in the sense of mixed Hodge theory.  (We conjecture this always
holds, cf. \S \ref{subsec-conjecture1}.)  The description in \S
\ref{thm-cohomology} is made possible by using \cite{GKM1} \S 6.3 (the lemma
of Chang and Skjelbred) which describes the equivariant homology of a
$T(\mathbb C)$-variety in terms of the $0$ and $1$ dimensional orbits of
$T(\mathbb C).$  The computation is reduced to the case of the equivariant
homology of affine Springer fibers for the group $\text{SL}(2).$  Among the
various ways to compute the equivariant homology in this special case of
$\text{SL}(2),$ we have chosen to describe one which is conceptually simple,
but which is computationally rather cumbersome.   It uses the existence of a
slightly larger torus $\widetilde{T} (\mathbb C)$ which acts on this Springer
fiber.  To facilitate the exposition we have extracted the combinatorial
details and placed them in \S \ref{sec-combinatorial}.

Suppose $G$ is adjoint, $s \in \widehat{T}$ is an element of the dual torus
and $(H,s)$ is endoscopic data for $G.$  Then $H$ and $G$ share the same torus
$T$ so $\gamma$ corresponds to a regular semisimple element $\gamma_H \in H,$
and there is an associated affine Springer fiber $X^H_{\gamma_H}.$  Motivated
by the fundamental lemma (\cite{Langlands} \S III.1), one predicts that there
should be a close relation between the homology $H_*(X_{\gamma};\bC)$ and
$H_*(X^H_{\gamma_H};\bC).$  However from a geometrical point of view, any
relationship between these homology groups comes as a surprise, since the
group $H$ may not be a subgroup of $G.$  In fact, the varieties $X_{\gamma}$
and $X^H_{\gamma_H}$ have very little to do with each other: there does not
appear to be any map or correspondence between them, and even their dimensions
are different, in general.  Moreover, the relation between these homology
groups turns out to be rather subtle:  they are not simply isomorphic, but
rather, they become isomorphic only after a degree shift and a certain
localization.

The element $s$ gives rise to a homomorphism $\chi_*(T) \to \mathbb C^{\times}
$ (defined on the cocharacter group of $T$) so that the set of co-roots for
$H$ is given by
\[ \Phi^{\vee}(H,T) = \left\{ \alpha^{\vee} \in \Phi^{\vee}(G,T):
s(\alpha^{\vee}) = 1 \right\}.\]
The same cocharacter group $\chi_*(T)$ acts freely on $X_{\gamma}$.  Let $J$
be the multiplicative subset of the group algebra $\mathbb C[\chi_*(T)]$ which
is generated by the elements $(1-\alpha^{\vee})$ where $\alpha^{\vee} \in
\Phi^{\vee}(G,T)$ and $s(\alpha^{\vee}) \ne 1.$  For any $\mathbb C[\chi_*(T)]
$ module $V,$ let $V_J = J^{-1}V$ be its localization with respect to $J.$ Let
$r = r(\gamma)$ be the non-negative integer defined in \S
\ref{subsec-endoscopic1}. The second principal result in this paper is the
following complex analog to the fundamental lemma in the unramified case (a
more precise version of which will be stated in Theorem \ref{thm-mainA}).

\begin{thm}
Suppose the homology of $X_{\gamma}$ is pure and the homology of
$X^H_{\gamma_H}$ is pure.  Then for each $i$ there is a homomorphism
\begin{equation}\label{eqn-themap}
H^{T(\mathbb C)}_i(X_{\gamma}, \mathbb C) \to H^{T(\mathbb C)}_{i-2r}
(X^H_{\gamma_H}, \mathbb C) \end{equation}
which becomes an isomorphism
\begin{equation}\label{eqn-FL0}
H_i^{T(\mathbb C)}(X_{\gamma}; \mathbb C)_J \cong H_{i-2r}^{T(\mathbb C)}
(X^H_{\gamma_H}; \mathbb C)_J
  \end{equation}
after localizing with respect to $J.$\end{thm}
The homomorphism (\ref{eqn-themap}) is compatible with a number of algebras
which act on these homology groups.  First, the equivariant homology is a
module over the equivariant cohomology $\mathcal D = H^*_{T(\mathbb C)}
(\text{pt})$ of a point (\S \ref{sec-equivariant}).  Second, the group algebra
$\mathbb C[\chi_*(T)]$ acts.  Finally, a certain group
$\widetilde{\mathcal W}^{G,H}_{\gamma}$
of automorphisms of the situation also acts (cf. \S \ref{subsec-Weylaction} and
\S \ref{subsec-endoscopic1}), 	and the homomorphism
(\ref{eqn-themap}) transforms by a certain
character $\eta$ under  $\widetilde{\mathcal W}^{G,H}_{\gamma}$.

Using the action of $\mathcal D$, it is possible to recover the ordinary
homology of $X_{\gamma}$ from its equivariant homology, so (\ref{eqn-themap})
gives a homomorphism
\[ H_i(X_{\gamma}; \mathbb C) \to H_{i-2r}(X^H_{\gamma_H}; \mathbb C)\]
of the ordinary homology, which also becomes an isomorphism after localizing
with respect to $J.$ In \S \ref{sec-endoscopic} the action of the group
algebra $\mathbb C[\chi_*(T)]$ is used in order to (partially) describe the
homology of the quotient $\chi_*(T) \backslash X_{\gamma}.$ Finally, the
action of the group $\widetilde{\mathcal W}^{G,H}_{\gamma}$ provides
a complex analog for the
action of the Frobenius automorphism.

These extra ingredients are related to the fundamental lemma for 
function fields. In Theorem
\ref{orb-thm}, the $\kappa$-orbital integral which occurs in the
fundamental lemma is shown to equal the trace of the Frobenius automorphism on
the \'etale cohomology of the quotient $\chi_*(T) \backslash X_{\gamma}
(\overline{k});$ where now $X_{\gamma} (\overline{k})$ is the finite
characteristic analog of the complex algebraic variety $X_{\gamma}.$  Next,
the analog of Theorems \ref{thm-cohomology} and \ref{thm-mainA} need to be
established in \'etale homology, a task which we have not fully carried out
here.  In \S \ref{subsec-remarksFL} we indicate how these results are related
to the fundamental lemma in the unramified case.

We would also like to draw attention to the recent preprint
\cite{LaumonFL} in which G. Laumon uses the methods and
techniques of the present paper together with a deformation
argument from \cite{Jacob} to prove (under
the same purity assumption) the ``geometric fundamental
lemma" (that is, the \'etale cohomology analog of Proposition
\ref{prop-E2}) for unitary groups (and for arbitrary regular
semi-simple elements $\gamma$).

In \S \ref{sec-flagSL2} and \S \ref{sec-flag} we list the changes which are
needed in order to establish similar results for Springer fibers $Y_{\gamma}$
in the affine flag manifold $Y.$  The homology of a Springer fiber in the
affine flag manifold carries the additional structure of a (right) action of
the affine Weyl group, the {\it Springer action}.  It was constructed by
Lusztig \cite{lusztig} and (using a statement from \cite{KL0} whose proof
does not appear in the literature) by Sage \cite{Sage2}.
In \S \ref{subsec-Springeraction} this action is explicitly described using our
formula (Theorem \ref{thm-affineflaghomology}) for the homology of the
Springer fiber $Y_{\gamma}$ (under the assumption that this homology is pure).
In Theorem \ref{thm-FLforflags} we state the main consequence:
If $(H,s)$ is endoscopic data for $G,$  if $H_*(Y_{\gamma}; \bC)$ is pure and
$H_*(Y^H_{\gamma_H};\bC)$ is pure, then for each $i$ we obtain a homomorphism
\[ H_i^{T(\mathbb C)}(Y_{\gamma};\bC) \to H_{i-2r}^{T(\mathbb C)}
(Y^H_{\gamma_H};\bC)\otimes_{\bC[\widetilde W^H]}\bC[\widetilde W]
\]
of $\mathcal D$ modules, and a homomorphism on ordinary homology,
\[
H_i(Y_{\gamma}; \mathbb C) \to H_{i-2r}(Y^H_{\gamma_H};\mathbb C)
\otimes_{\bC[\widetilde W^H]}\bC[\widetilde W],
\]
each of which is equivariant with respect to the right action of
$\widetilde{W}$ and transforms by the character $\eta$ under the action of
$\widetilde{\mathcal W}^{G,H}_{\gamma}$.
Each of these homomorphisms becomes an isomorphism after localizing with
respect to $J.$

We wish to thank R. Hain and David Nadler for useful conversations.
The authors are grateful to the University of Chicago, the
Universit\'e de Paris VII, the \'Ecole Normale Sup\'erieure, the Institut
Henri Poincar\'e, and the Institute for Advanced Study for their support
during various stages of this research.

\section{Purity for {\it ind}-varieties}\label{sec-purity}
\subsection{}
Let $X$ be a complex projective algebraic variety with (increasing) weight
filtration $W$ of the cohomology $H^*=H^*(X;\mathbb C)$ (see \cite{Deligne}).
The cohomology of $X$ is {\it pure} if
\[\text{Gr}^m W(H^i) = W^m(H^i)/W^{m-1}(H^i) = 0\ \text{ for } i \ne m.\]
We obtain an increasing weight filtration on the homology
$H_*=H_*(X; \mathbb C):$
\[
W^m(H_i) = \left( H^i/W^{-1-m}(H^i) \right)^*
\]
by duality.  If $f:X \to Y$ is a morphism of projective varieties, the induced
mapping $f_*:H_*(X) \to H_*(Y)$ is strictly compatible with the weight
filtration.

Let $X$ be a complex projective {\it ind}-variety (\cite{Kumar},
\cite{Shafarevich}), that is, a set with a filtration $X_0 \subset X_1 \subset
\ldots$ by complex projective varieties such that $X = \bigcup_{n=1}^{\infty}
X_n$ and such that each $X_n \to X_{n+1}$ is a closed immersion.  A set $Y
\subset X$ is closed in the limit topology for $X$ iff $Y \cap X_n$ is closed
(in the classical topology) for each $n.$  Then $H_*(X) =
\underset{\rightarrow} {\lim} H_*(X_n)$ inherits a limit weight filtration $W$
from the weight filtrations on $H_*(X_n).$  If $X = \bigcup_{n=1}^{\infty}
Y_n$ is an equivalent {\it ind}-variety structure on $X$ (meaning that the
identity mapping $X \to X$ is a morphism of {\it ind}-varieties) then the
resulting weight filtration $W'$ on $H_*(X)$ agrees with $W.$  We say the
homology of $X$ is pure if $\text{Gr}^{-m}W(H_i(X)) =0$ for $m \ne i.$

Let $A= (\mathbb C^{\times})^m$ be a complex torus which acts on a projective
{\it ind}-variety $X$ such that each $X_n$ is invariant under $A$ and so that
the action of $A$ on $X_n$ is algebraic and is compatible with the immersion
$X_n \to X_{n+1}.$ In this case we say that $A$ acts algebraically on $X.$
The classifying space $BA$ may be taken to be the {\it ind}-variety $(\mathbf
P^1)^m \subset (\mathbf P^2)^m \ldots,$ which has pure homology.  The
Cartan-Leray spectral sequence for the $A$-equivariant homology of $X_n$ is
the spectral sequence for the fibration $X_n \times_A EA \to BA$ with $E^2
\cong H_*(BA) \otimes H_*(X_n).$  It inherits a weight filtration such that
the differentials are strictly compatible with the weight.  It follows that
the Cartan-Leray spectral sequence for the $A$-equivariant homology of $X$
inherits a weight filtration such that the differentials are strictly
compatible with the weight.  The following fact will be needed for Corollary
\ref{cor-ordinaryH}.

\begin{lem}\label{lem-pureind}
Suppose $X$ is a complex projective {\it ind}-variety with an action of an
algebraic torus $A.$ Suppose the homology of $X$ is pure.  Then the homology
of $X$ is equivariantly formal \cite{GKM1} \S 1, that is, the Cartan-Leray
spectral sequence for the $A$-equivariant homology of $X$ collapses at $E^2$
and gives rise to an isomorphism
\[ H^A_*(X;\mathbb C) \cong H_*(X;\mathbb C)\otimes_{\mathbb C}
H^A_*(\text{pt}; \mathbb C).\]
\end{lem}

\subsection{Purity of $B\mathbb Z^n$}
Define an action of the additive group $\mathbb Z$ on the
contractible {\it ind}-variety
\begin{equation*}
E\mathbb Z = (\mathbb Z \times \mathbb C)/(m,0)\sim (m+1,1)
\end{equation*}
(which is a countable union of affine lines joined each to the
next at a single point) by $m \cdot (r,t) = (r+2m,t)$ (where $m,
r \in \mathbb Z$ and $t\in \mathbb C$).   The quotient $B\mathbb
Z = E\mathbb Z / \mathbb Z$ consists of two copies $V_0, V_1$ of
$\mathbb C,$ joined at two points:  $0\in V_0$ is identified with
$1\in V_1$ and vice versa.  The inclusion of the unit circle $T^1
\to B\mathbb Z$ given by
\begin{equation*}
e^{i\theta} \mapsto \begin{cases}
\frac{\theta}{\pi}\in V_0 &\text{ if }\ 0 \le \theta \le \pi\\
\frac{\theta}{\pi}-1\in V_1 &\text{ if }\ \pi \le \theta \le 2\pi
\end{cases}\end{equation*}
is a homotopy equivalence.  Let $\Lambda = \mathbb Z^n,$ and set
$E\Lambda = (E\mathbb Z)^n$ and $B\Lambda = (B\mathbb Z)^n.$  We
obtain a homotopy equivalence $T^n = (T^1)^n \to B\Lambda.$

Suppose $X$ is an {\it ind}-variety on which $\Lambda$ acts
freely by morphisms.  Let $E(X) = E\Lambda\times_{\Lambda} X$ be
the Borel construction, together with its projections
\begin{equation*} \begin{CD}
\Lambda \backslash X @<{p}<< E(X) @>{\pi}>> B\Lambda \end{CD}.
\end{equation*}

Let $s: \Lambda \to \mathbf{GL}(1,\mathbb C)$  be a 1-dimensional
representation of finite order, and write  $\mathbb C_s$ for the
representation space.  Then $s$ determines local systems
$\mathcal L_s = \mathbb C_s \times_{\Lambda} X$ on $\Lambda
\backslash X$ and  $\mathcal L^B_s = E\Lambda \times_{\Lambda}
\mathbb C_s$ on $B\Lambda,$ and there is a canonical isomorphism
of local systems, $\pi^*\mathcal L^B_s \cong p^*\mathcal L_s.$
Denote this local system on $E(X) $ by $\mathcal L^E_s.$  The
homology Leray spectral sequence, with coefficients in $\mathcal
L^E_s,$ for the fibration $\pi:E(X) \to B\Lambda$ has
\begin{equation}\label{eqn-E2a}
E^2_{pq} = H_p(\Lambda; H_q(X; \mathbb C_s))
= H_p(B\Lambda; \mathcal L^B_s\otimes \mathbf{H_q}) \implies
H_{p+q}(\Lambda\backslash X; \mathcal L_s) \end{equation}
where $\mathbf {H_q}$ denotes the local system on $B\Lambda$
corresponding to $H_q(X;\mathbb C).$  The following fact will be
needed for Proposition \ref{prop-E2}.
\begin{lem}\label{lem-collapse}
  Suppose the homology of $X$ is pure.  Then the spectral
sequence \parens{{\rm\ref{eqn-E2a}}} collapses at $E^2.$  Hence
there is an isomorphism
\begin{equation}\label{eqn-directsum}
H_m(\Lambda\backslash X; \mathcal L_s) \cong \bigoplus_{p+q=m}
H_p(\Lambda;H_q(X; \mathbb C_s)).  \end{equation}\end{lem}

\subsection{Proof}

First consider the case when $n=1$ and when
$s$ is the trivial representation.  The homology $H_*(B\mathbb Z;
\mathbf{H_q})$ may be computed from the Mayer Vietoris
sequence for the covering of $B\mathbb Z$ by the two (Zariski)
closed sets $V_0$ and $V_1.$  It is easy to see that this sequence
gives a short exact sequence (of mixed Hodge structures)
\begin{equation*}\begin{CD}
0 @>>> H_1(B\mathbb Z; \mathbf{H_q}) @>>> H_q @>>{\mu-1}> H_q @>>>
H_0(B\mathbb Z; \mathbf{H_q}) @>>> 0
\end{CD}\end{equation*}
where $H_q = H_q(X; \mathbb C),$ and where $\mu$ is the
homomorphism induced on homology from the action of $1 \in
\mathbb Z.$

Next consider the case of general $n$ but still with $s$ the
trivial representation.  It follows by induction that the
homology $H_*(B\mathbb Z^n; \mathbf{H_q})$ is the homology of the
Koszul complex $K_{\bullet} = \textstyle\bigwedge^*(\mathbb C^n)
\otimes H_q$ with differential \begin{equation*}
\partial(e_{i_1}\wedge\ldots \wedge e_{i_r}\otimes h) = \sum_{j=
1}^r (-1)^{j+1}e_{i_1}\wedge \ldots\wedge \widehat{e}_{i_j}\wedge
\ldots
\wedge e_{i_r} \otimes (\mu_{i_j}(h)-h)
\end{equation*}
where $\mu_1,\ldots,\mu_n$ are the homomorphisms on $H_q$
determined by the action of the standard $n$ basis elements $e_1,
\ldots, e_n.$  Each $\textstyle\bigwedge^k(\mathbb C^n) \otimes
H_q$ is pure of weight $-q$ since it arises as a direct sum of
copies of $H_q.$  So each homology group $H_p(B\Lambda; \mathbf{H_q})$
is pure of weight $-q.$   Since the differentials in the spectral
sequence (\ref{eqn-E2a}) are strictly compatible with the weight
filtration, this implies the spectral sequence collapses at
$E^2.$

Now consider the case of a nontrivial representation $s$
of finite order.  Such a representation factors through the
(finite) quotient $\Lambda/N\Lambda$ for some positive integer
$N.$  Let $\phi:\widetilde{B}= E\Lambda / N \Lambda \to B\Lambda$
be the resulting (finite) cover, and let $\widetilde{E}$ be the
fiber product in the following diagram,
\begin{equation*}\begin{CD}
\widetilde{E} @>{\Phi}>> E\Lambda\times_{\Lambda}X @>>{p}>
\Lambda \backslash X \\
@V{\tilde{\pi}}VV @VV{\pi}V && \\
\widetilde{B} @>>{\phi}> B\Lambda
\end{CD}\end{equation*}
The preceding argument may be used to see that the Leray spectral
sequence for the mapping $\tilde\pi$ (with constant coefficients)
collapses at $E^2.$

The covering $\phi$ is a principal homogeneous space for the
group $\Lambda/N\Lambda,$ which acts by deck transformations on
$\widetilde{B}.$  So it acts on the local system $\phi_*(\mathbb
C)$ as the regular representation, and determines decompositions
\begin{equation}\label{eqn-ch}
\phi_*(\mathbb C) \cong \underset{\theta}{\textstyle\bigoplus}\
\mathcal L^B_{\theta}\ \text{ and }\
\Phi_*(\mathbb C) \cong \underset{\theta}{\textstyle\bigoplus}\
\mathcal L^E_{\theta} \end{equation}
into one dimensional local systems corresponding to the distinct
characters $\theta: \Lambda/N\Lambda \to \mathbb C^{\times}.$
It follows (see below) that there is a natural isomorphism
\begin{equation}\label{eqn-hd}
H_*(\widetilde{B}; \mathbf{H_q^{\widetilde{B}}}) \cong
\underset{\theta}{\textstyle\bigoplus}\ H_*(B\Lambda;
\mathbf{H_q}\otimes \mathcal L_{\theta}) \end{equation}
where $\mathbf{H_q^{\widetilde{B}}}$ (resp.~$\mathbf{H_q}$) is
the local system on $\widetilde{B}$ (resp.~on $B\Lambda$)
corresponding to $H_q(X;\mathbb C).$  In fact the whole spectral
sequence for $\tilde\pi$ decomposes under $\Lambda/N\Lambda$ into
a direct sum (over distinct characters $\theta$ of
$\Lambda/N\Lambda$) of spectral sequences for $\pi$ with
coefficients in $\mathcal L_{\theta}.$  It follows that each of
these constituent spectral sequences collapses, one of which
corresponds to the character $s.$

There are several ways to verify equation (\ref{eqn-hd}).  One
way is to use the homotopy equivalence $T^n \to B\Lambda,$
Poincar\'e duality, and (\ref{eqn-ch}).  However it may also be
verified directly.  Let $\mathbf {D_{B\Lambda}}$
(resp.~$\mathbf{D_{\widetilde{B}}}$) be the dualizing complex on
$B\Lambda$ (resp.~on $\widetilde{B}$).  Since $\phi$ is a finite
covering space, $\phi_* \cong R\phi_* \cong R\phi_{!}$ and
$\phi^*(\mathbf{D_{B\Lambda}}) \cong \mathbf D_{\widetilde{B}}.$
Since the diagram above is a fiber square,
\begin{equation*}
R\phi_*(\mathbf H_q^{\widetilde{B}}) \cong
R\phi_*(\phi^*\mathbf{H_q}\otimes\mathbb C) \cong \mathbf {H_q}
\otimes \left( \underset{\theta}{\textstyle\bigoplus}\ \mathcal
L_{\theta}^B \right).
\end{equation*}
The homology $H_*(\widetilde{B}; \mathbf{H_q^{\widetilde{B}}})$
is the cohomology (with compact supports) of the sheaf
\begin{equation*}
\mathbf{D_{\widetilde{B}}}\otimes \mathbf{H_q^{\widetilde{B}}}
\cong
\phi^*\mathbf{D_B}\otimes \mathbf{H_q^{\widetilde{B}}}.
\end{equation*}
Pushing this sheaf forward under $\phi$ gives
\begin{equation*}
\mathbf{D_B} \otimes R\phi_*(\mathbf{H_q^{\widetilde{B}}})
\cong \mathbf {D_B} \otimes \mathbf{H_q} \otimes \left(
\underset{\theta}{\textstyle\bigoplus}\ \mathcal L_{\theta}
^B\right) \end{equation*}
whose compactly supported cohomology is $\underset{\theta}
{\textstyle\bigoplus}\ H_*(B\Lambda; \mathbf{H_q}\otimes  
\mathcal L_{\theta}^B)$ as claimed. \qed

\subsection{Remarks}
Presumably the theory of weights can be extended to include
(co)homology groups such as $H_k(B\Lambda; \mathbf{H_q}\otimes
\mathcal L_s),$ in which case equation (\ref{eqn-hd}) implies
that this group is pure of weight $-q.$   We remark that if
the homology of $X$ is both pure and {\it Tate} (meaning that
it lives only in even degrees and that $Gr^{-2q}H_*(X)$ is all
of Hodge type $(q,q)$), then the decomposition
(\ref{eqn-directsum}) is canonical because in this
case, the Hodge filtration and the weight filtration are opposed
(\cite{HodgeII} Prop. 1.2.5).  It seems likely that Lemma
\ref{lem-collapse} remains valid even when the local system $s$
fails to have finite order.

\section{Symmetric algebras}\label{sec-symmetricalgebras}
\subsection{}
Consider the polynomial algebra
\[\bS = \mathbb Q[x_1,\ldots,x_n] = \textstyle{\bigoplus_{k=0}^{\infty}}
\bS_k\]
in $n$ variables, graded by degree (with finite dimensional graded pieces);
and its dual algebra of differential operators
\[\mathcal D = \mathbb Q[\partial_1,\ldots,\partial_n] =
\textstyle{\bigoplus_{k=0}^{\infty}}
\mathcal D_k\]
with $\partial_i x_j = \delta_{ij}.$  Then $\mathcal D_k \bS_j \subset
\bS_{j-k}.$  The natural pairing
\begin{equation}\label{eqn-pairing} \begin{CD}
\mathcal D \otimes \bS @>{\langle \cdot ,\cdot\rangle}>> \mathbb Q\quad
\end{CD}\end{equation}
given by $\langle \partial, P \rangle = (\partial P)(0)$ satisfies
\begin{equation}\label{eqn-adjoint}
\langle \partial \partial',P \rangle = \langle \partial, \partial'P \rangle.
\end{equation}
It restricts to a nondegenerate (finite dimensional) bilinear form $\mathcal
D_k \otimes \bS_k \to \mathbb Q,$ which vanishes on $\mathcal D_k \otimes
\bS_j$ for $j \ne k.$ If $I = \bigoplus_{k=0}^{\infty} I_k \subset \mathcal D$
is a homogeneous ideal, define
\[ \bS\{I\} = \left\{ s\in\bS:\ \partial s = 0 \text{ for all } \partial \in I
\right\}
\]
to be the subspace annihilated by $I.$  It follows from (\ref{eqn-adjoint})
that
\begin{equation}\label{eqn-perp}
\bS\{I\} = I^{\perp}=\textstyle{\bigoplus_{k=0}^{\infty}} I_k^{\perp}
\end{equation}
where
\begin{equation*}
I_k^{\perp} = \left\{ s \in \bS_k:\ \langle \partial, s \rangle = 0 \text{ for
all } \partial \in I_k \right\}.  \end{equation*}
So the pairing (\ref{eqn-pairing}) passes to a nondegenerate pairing
$(\mathcal D_k/I_k) \otimes I_k^{\perp} \to \mathbb Q.$ If $\partial \in
\mathcal D$ set $\bS\{\partial\} = \bS\{(\partial)\}.$

\begin{lem}\label{lem-relativelyprime}
Suppose $\partial, \partial' \in \mathcal D$ are relatively prime homogeneous
elements of $\mathcal D.$  Then
\begin{equation}\label{eqn-sum}
\bS\{\partial \partial'\} = \bS\{\partial\} + \bS\{\partial'\}
\end{equation}
and
\begin{equation}\label{eqn-prime} \begin{CD}
\bS\{\partial\} @>{ \partial'\cdot}>> \bS\{\partial\}
\end{CD} \end{equation}
is surjective and lowers degree by $\deg(\partial').$ \end{lem}

\subsection{Proof}
If $I,J \subset \mathcal D$ are homogeneous ideals, then it follows from
(\ref{eqn-perp}) that
\begin{align*}
\bS\{I+J\} &= \bS\{I\} \cap \bS\{J\} \\
\bS\{I\cap J\} &= \bS\{I\} + \bS\{J\}.
\end{align*}
Taking $I= (\partial)$ and $J = (\partial')$ proves (\ref{eqn-sum}).  Since
$\partial$ and $\partial'$ are relatively prime, the mapping
\begin{equation*}\begin{CD}
\mathcal D/(\partial) @>{ \partial'\cdot}>> \mathcal D/(\partial)
\end{CD}\end{equation*}
is an injective homomorphism of graded algebras, which raises degree by
$\deg(\partial').$  It follows by duality (\ref{eqn-pairing}) that
(\ref{eqn-prime}) is surjective.
\qed

\section{Equivariant homology}\label{sec-equivariant}
\subsection{}\label{subsec-equivariant1}
  Let $A=A(\mathbb C)$ be a complex torus with Lie algebra $\mathfrak a.$ We
will often use without mention the canonical isomorphism $\chi^*(A) \otimes
\mathbb C \to \mathfrak a^*.$  Let
\[
\mathcal D(\mathfrak a) = \text{Sym}(\mathfrak a^*) = \bigoplus_{d=0}^{\infty}
\text{Sym}^d(\mathfrak a^*)\]
be the symmetric algebra of linear differential operators (with constant
complex coefficients) on $\mathfrak a^*,$ graded by degree.  Let
\[
\bS(\mathfrak a) = \text{Sym}(\mathfrak a) = \bigoplus_{d=0}^{\infty}
\text{Sym} ^d(\mathfrak a^*) ^* \]
be the dual symmetric algebra of complex valued polynomial functions on
$\mathfrak a^*.$

Throughout this paper, cohomology and homology will be taken with complex
coefficients (unless otherwise stated), however it is possible to use rational
coefficients, (by replacing $\mathfrak a^*$ with $\chi^*(A) \otimes \mathbb
Q,$ for example).  Let $BA$ be the classifying space for $A.$  If $\varphi: A
\to \mathbb C^{\times}$ is a character  let $c_1(\varphi)$ be the first Chern
class of the resulting line bundle on $BA.$ Then $c_1$ extends to an
isomorphism
\[
\mathcal D(\mathfrak a) \to H^*(BA) = H^*_A(\text{pt})\]
(the Chern-Weil isomorphism) of graded algebras, which doubles degrees.

Multiplication $A \times A \to A$ determines an H-space structure $m: BA
\times BA \to BA$ which gives rise to an algebra structure on $H_*(BA) =
H_*^A(\text{pt}).$  The homology $H_*(BA)$ is also a module over $H^*(BA)$
under the cap product.   If $x \in H^2(BA)$ then the Hopf formula $m^*(x) =
1\otimes x + x \otimes 1$ implies that $x$ acts on $H_*(BA)$ as a derivation,
so $H^*(BA)$ acts on $H_*(BA)$ by differential operators.  It follows that the
mapping $\bS(\mathfrak a)  \to H_*(BA)$ (dual to the Chern-Weil isomorphism)
is an algebra isomorphism and is compatible with the cap product, in the sense
that the following diagram commutes:
\begin{equation*} \begin{CD}
\mathcal D_k \otimes \bS_j @>>> \bS_{j-k} \\
@VV{\cong}V @VV{\cong}V \\
H^{2k}(BA) \otimes H_{2j}(BA) @>{\cap}>> H_{2j-2k}(BA)
\end{CD} \end{equation*}

Suppose $A$ acts on a complex projective algebraic variety $X$ whose
cohomology is pure.  Then the cohomology of $X$ is equivariantly formal and
the Cartan-Leray spectral sequence for the equivariant cohomology of $X$
collapses at $E_2.$  Let $T\subset A$ be a (possibly trivial) subtorus and let
$I = \ker(H^*(BA) \to H^*(BT))$ be the resulting homogeneous ideal.  Then the
$T$-equivariant cohomology of $X$ may be recovered from the $A$-equivariant
cohomology by
\[ H^*_T(X) = H^*_A(X) \otimes_{H^*_A(\text{pt})} H^*_T(\text{pt}).\]
The equivariant homology $H_*^A(X)$ is a module over $\mathcal D =
H^*_A(\text{pt})$ and by duality, the $T$-equivariant homology of $X,$
\begin{equation}\label{eqn-annihilator}
  H_*^T(X) = H_*^A(X)\left\{ I \right\} \end{equation}
is the submodule of $H_*^A(X)$ consisting of those elements which are
annihilated by the homogeneous ideal $I.$

\subsection{}\label{subsec-diagram}
Let $\varphi:A \to \mathbb C^{\times}$ be a nontrivial character with
corresponding differential operator $\partial_{\varphi} \in \mathcal
D(\mathfrak a).$  Let $M = \ker(\varphi)$ with Lie algebra $\mathfrak m =
\ker(\partial_{\varphi})\subset \mathfrak a^*.$ Let $j^*: \mathcal S(\mathfrak
m) \to \bS(\mathfrak a)$ be the homomorphism induced by inclusion $\mathfrak m
\subset \mathfrak a.$  Then the image of $j^*$ consists of the polynomial
functions $\mathcal S(\mathfrak a)\left\{ \partial_{\varphi} \right\}$ on
$\mathfrak a^*$ which are annihilated by $\partial_{\varphi}.$

Such a nontrivial character $\varphi$ extends to an action of $A$ on $\mathbb
C \mathbb P^1 = \mathbb C^{\times} \cup \left\{ 0 \right\} \cup \left\{ \infty
\right\}.$  As in \cite{GKM1} \S 7.1, 7.2, the long exact sequence for the
pair $(\mathbb C \mathbb P^1, \left\{ 0 \right\} \cup \left\{ \infty \right\})
$ breaks into short exact sequences
\begin{equation*}
\begin{CD}
0 @>>> H_*^A(\mathbb C \mathbb P^1, \left\{ 0 \right\} \cup \left\{ \infty
\right\} )
@>{\partial}>> H_*^A(\left\{ 0 \right\}) \oplus H_*^T(\left\{ \infty \right\})
@>>> H_*^A(\mathbb C \mathbb P^1) @>>> 0  \\
&&||&&|| \\
&& \bS(\mathfrak m) @>>{\beta}> \bS(\mathfrak a) \oplus \bS(\mathfrak a)
\end{CD} \end{equation*}
where $\beta(f) = (j^*(f), -j^*(f)).$

\subsection{}\label{subsec-equivariant2}
Now suppose that a complex algebraic torus $A$ acts algebraically on a complex
projective variety $Y$ whose cohomology is pure.  Let $Y_0 \subset Y$ be the
set of fixed points and let $Y_1\subset Y$ be the union of the 0 and 1
dimensional orbits of $A$ in $Y.$  The lemma of Chang and Skjelbred \cite{CS}
(cf. \cite{GKM1} \S 6.3), when translated into homology, says that the
following sequence is exact:
\begin{equation}\label{eqn-localization}
\begin{CD}
H^A_*(Y_1, Y_0) @>{\partial}>> H^A_*(Y_0) @>>> H^A_*(Y) @>>> 0.
\end{CD} \end{equation}
Suppose, moreover, that there are finitely many fixed points $Y_0 = \left\{
y_1,y_2,\ldots, y_r\right\}$ and finitely many 1-dimensional orbits $\left\{
E_1, E_2, \ldots, E_d \right\}.$ The closure $\overline{E}_i$ of each
1-dimensional orbit is obtained by adding two fixed points $\partial E_i =
\left\{ y_{i_a} \right\} \cup \left\{ y_{i_b} \right\}.$  Let $\mathfrak m_i
\subset \mathfrak a$ be the Lie algebra of the stabilizer of any point in
$E_i$ and let $j_i^*: \bS(\mathfrak m_i) \to \bS(\mathfrak a)$ be the
resulting homomorphism.  Using \S \ref{subsec-diagram}, the sequence
(\ref{eqn-localization}) becomes
\begin{equation}\label{eqn-homologysequence} \begin{CD}
\bigoplus_{i=1}^d \bS(\mathfrak m_i) @>>{\beta}>
\bigoplus_{k=1}^r \bS(\mathfrak a) @>>> H_*^A(Y) @>>> 0
\end{CD} \end{equation}
where $\beta = \sum_{i=1}^d \beta_i$ with
\[ \beta_i(f_i) = \left( 0,\ldots, j_i^*(f_i), \ldots, -j_i^*(f_i), \ldots, 0
\right).\]
for any $f_i \in \bS(\mathfrak m_i).$  (Here, the two nonzero entries occur in
the positions corresponding to $y_{i_a}$ and $y_{i_b}$ respectively.)

\section{Affine Springer fibers}\label{sec-affinespringer}
\subsection{Notation}\label{subsec-preliminaries}
Let $\mathbb G_m$ denote the multiplicative group.  If $k$ is a field, if
$\overline{k}$ is an algebraic closure, and if $T$ is an algebraic torus
defined over $k$, let $\chi_*(T) = \text{Hom}(\mathbb G_m,T) \cong \text{Hom}
(\overline{k}^{\times}, T(\overline{k}))$ denote the group of cocharacters of
$T$ and let $\chi^*(T) = \text{Hom}(T, \mathbb G_m) \cong \text{Hom}
(T(\overline{k}), \overline{k}^{\times})$ denote the group of characters of
$T$. The dual torus is $\widehat{T} = \text{Hom}(\chi_*(T), \mathbb G_m).$
If $G$ is a connected reductive complex algebraic group, if $T\subset G$ is a
maximal torus in $G$ and if $(\chi^*(T),\Phi, \chi_*(T), \Phi^{\vee})$ is the
resulting root datum for $G$ (where $\Phi$ and $\Phi^{\vee}$ are the roots and
co-roots of $T$ in $G$ respectively), then we denote by $\widehat{G}$ the
``dual'' connected reductive complex algebraic group corresponding to the root
datum $(\chi_*(T), \Phi^{\vee}, \chi^*(T), \Phi).$

Throughout this paper (except in \S \ref{sec-orbitalintegral}) we let $F =
\mathbb C((\varepsilon)) $ be the field of formal Laurent series over $\mathbb
C$ and denote by $\mathfrak o = \mathbb C[[\varepsilon]] $ its ring of
integers, the formal power series over $\mathbb C.$  The valuation
$\text{val}:F \to \mathbb Z$ takes integer values.  Let $T$ be an algebraic
torus defined over $\mathbb C.$  The following sequence is exact:
\begin{equation*} \begin{CD}
  1 @>>> T(\mathfrak o) @>>> T(F) @>{\text{val}}>> \chi_*(T) @>>> 1.
\end{CD}\end{equation*}
Here, $\val{}$ is defined by the property that $\alpha(\text{val}(\ell)) =
\text{val} (\alpha(\ell))$ for all $\alpha \in \chi^*(T)$ and all $\ell \in
T(F).$  (On the left side of this equality, $\alpha$ is viewed as a
homomorphism $\chi_*(T) \to \mathbb Z.$) The choice of uniformizing parameter
$\varepsilon$ determines a splitting $\chi_*(T) \to T(F)$ by $\beta \mapsto
\beta(\varepsilon),$ whose image $\Lambda$ will be called the lattice of
translations.  It is a free abelian group of rank equal to the dimension of
$T.$

\subsection{} In this section we recall some definitions and results from
\cite{KL}. Let $G$ be a connected reductive algebraic group over $\mathbb C,$
with Lie algebra $\mathfrak g$.  Set $\mathfrak g(F) = \mathfrak g
\otimes_{\mathbb C} F$ and $\mathfrak g({\mathfrak o}) = \mathfrak g \otimes
\mathfrak o.$  Denote by $\text{Ad}$ the adjoint action of $G_F$ on $\mathfrak
g_F.$  We will often write $K$ for the group $G({\mathfrak o})$ of $\mathfrak
o$-points of $G$.

The {\it loop Grassmannian} (or affine Grassmannian) is the quotient $X = G(F)
/K.$ The affine Grassmannian $X$ is an {\it ind}-algebraic variety: it is an
increasing union \cite{KL}
\[ X_0 \subset X_1 \subset X_2 \subset \ldots \subset X\]  of finite
dimensional complex projective algebraic varieties.  If $H \subset G$ is a
connected reductive algebraic subgroup then the inclusion of $H$ into $G$
induces an injection $H(K)/H(\mathfrak o) \hookrightarrow X$ of the loop
Grassmannian for $H$ into the loop Grassmannian for $G.$

Each element $\gamma \in \mathfrak g(F)$ gives rise to a ``vector field'' on
$X$ whose fixed point set
\[ X_{\gamma} = \left\{ xK \in G(F)/K :\ \text{Ad}(x^{-1})(\gamma)
\in\mathfrak g({\mathfrak o})\right\} \]
is called an affine Springer fiber.  We say the element $\gamma$ is {\it
compact} if $X_{\gamma} \ne \phi$.  In \cite{KL} it is proven that a compact
element $\gamma \in \mathfrak g(F)$ is regular and semisimple iff $X_{\gamma}$
is finite dimensional, in which case $X_{\gamma}$ is an {\it ind}-subvariety
of $X.$  It is a union of (possibly) infinitely many irreducible components,
each of which is a complex projective variety.

\subsection{Conjecture}\label{subsec-conjecture1}  If $\gamma \in \mathfrak
g(F)$ is compact, regular, and semisimple then for all $i$ the homology group
$H_i(X_{\gamma};\mathbb C)$ is pure of weight $i.$

In \cite{GKM2} we prove this conjecture for elements $\gamma$ which have
``equal valuation''.

\subsection{Bruhat decomposition}
Let $T \subset G$ be a  maximal torus defined over $\mathbb C$ (hence split
over $F$), with its lattice of translations $\Lambda \subset T(F) \subset
G(F).$  The identification $\Lambda \cong T(F)/T(\mathfrak o)$ induces an
embedding $\phi:\Lambda \to G(F)/G(\mathfrak o) = X$ of the loop Grassmannian
for $T$ into the loop Grassmannian for $G.$  Fix a Borel subgroup $B \supset
T$ and let $I \subset G(F)$ be the corresponding Iwahori subgroup.  Then $G(F)
= I\Lambda K$ so $X$ decomposes into $I$ orbits or Bruhat cells,
\[ X = \coprod\limits_{\ell \in \Lambda} I \ell K.\]
We denote by $C_{\ell} = I \ell x_0 \subset X$ the orbit (cell) corresponding
to $\ell \in \Lambda.$

\subsection{Turning torus}\label{subsec-turning}
Each $\lambda \in \mathbb C^{\times}$ determines an automorphism
$\sigma_{\lambda}$ of the field $F = \mathbb C ((\varepsilon))$ by
$\sigma_{\lambda}(\varepsilon^m) = \lambda^m \varepsilon^m.$  We say that
$\lambda$ acts on $F$ by ``turning the loop'' and we refer to the image of
$\mathbb C^{\times} \to \text{Aut}(F)$ as the ``turning torus''.  The turning
torus preserves $\mathfrak o = \mathbb C [[\varepsilon]]$ and the fixed point
set of its action is
\[ \mathfrak o^{\mathbb C^{\times}} = F^{\mathbb C^{\times}} = \mathbb C \cdot
\varepsilon^0.\]
The resulting action of $\mathbb C^{\times}$ on $G(F)$ preserves $K =
G(\mathfrak o)$ and induces an action on $X = G(F)/K$ which commutes with the
action of $T(\mathbb C) \subset I.$ On each Bruhat cell $C_{\ell}$ the action
of the extended torus
\begin{equation}\label{eqn-extendedtorus}
  \widetilde{T}(\mathbb C) = T(\mathbb C) \times \mathbb C^{\times}
\end{equation}
may be described by
\begin{equation}\label{eqn-torusaction}
(t,\lambda)\cdot g\ell K = \lambda (tgt^{-1}) \ell K.  \end{equation}

\subsection{Affine roots}
   Let $\Phi = \Phi(G,T)$ be the root system for $G$, with positive roots
$\Phi^+$ determined by the choice of $B$, and with root space decomposition
\[ \mathfrak g(\mathbb C) = \mathfrak t(\mathbb C) \oplus \bigoplus_{\alpha
\in \Phi} \mathbb
C Y_{\alpha}\]
where $Y_{\alpha} \in \mathfrak g(\mathbb C)_{\alpha}$ are root vectors. Let
$\widetilde{\Phi} = \left\{ (\alpha, k):\  \alpha \in \Phi, k \in \mathbb Z
\right\}$ be the collection of affine roots.  Each affine root $(\alpha,k)$
may be considered to be a character of the extended torus $\widetilde{T}
(\mathbb C)$ (\ref{eqn-extendedtorus}) with $(\alpha,k)(t,\lambda) = \alpha(t)
\lambda^k.$  Then $\widetilde{T}(\mathbb C)$ acts on the affine root space
$\mathbb C \varepsilon^k Y_{\alpha}$ through this character, that is,
\begin{equation} \label{eqn-rootaction}
(t,\lambda)\cdot \varepsilon^k Y_{\alpha} = \alpha(t)
\lambda^k\varepsilon^kY_{\alpha}.\end{equation}

Let $C_0$ be the fundamental alcove determined by $I$ and let
$\widetilde{\Delta}$ be the set of  simple affine roots, so
\[
C_0 = \left\{ a \in \chi_*(T)\otimes_{\mathbb Z} \mathbb R :\
\alpha(a)+k > 0 \text{ for all } (\alpha,k) \in \widetilde{\Delta}  \right\}.
\]
There is a split short exact sequence
\[ 1 \to I_{+} \to I \to T(\mathbb C) \to 1\]
where $I_+$ denotes the unipotent radical of $I.$  Its Lie algebra is given by
\begin{equation}\label{eqn-rootspaces}
\mathfrak N(I) = \text{Lie}(I_+) = \prod_{n\in \mathbb Z}
\varepsilon^nT(\mathbb C) \oplus\prod_{(\alpha,k)}
  \mathbb C \varepsilon^k Y_{\alpha}  \end{equation}
where the product is taken over those affine roots $(\alpha,k) \in
\widetilde{\Phi}$ such that $\alpha(a)+k > 0$ for all $a \in C_0.$  To fix
notation we recall the following standard description of the Bruhat cells.

\begin{lem}\label{lem-Bruhatcell}  Let $x_0 = K \in X$ denote the basepoint.
Fix $\ell \in \Lambda.$   The exponential map determines a $\widetilde{T}
(\mathbb C)$-equivariant isomorphism between the Bruhat cell $C_{\ell} = I
\ell x_0 $ and the vector space
\begin{equation}\label{eqn-complement}
D_{\ell}=  \bigoplus\limits_{ (\alpha,k) } \mathbb C \varepsilon^k Y_{\alpha}
\end{equation}
where the sum is over those \parens{finitely many} affine roots $(\alpha,k)
\in \widetilde{\Phi}$ such that
\begin{equation*}
  \text{\rm val} (\alpha(\ell)) + k < 0 \text{ and }
\alpha(a) + k > 0  \text{ for all } a \in C_0
\end{equation*}
and where $\widetilde{T}(\mathbb C)$ acts on this vector space according to
\parens{\ref{eqn-rootaction}}.
\end{lem}

\subsection{Proof}  The subgroup $I_+$ acts transitively on $C_{\ell}.$  The
stabilizer $S_{\ell}$ of the point $\ell x_0$ is $I_+ \cap \ell G(\mathfrak o)
\ell^{-1}$  whose Lie algebra $\mathfrak s_{\ell}$ is the sum of the affine
root spaces $\mathbb C \varepsilon^kY_{\alpha}$ such that $\alpha(a)+k >0$
(for $a \in A$) and $\text{Ad}(\ell)(\varepsilon^k Y_{\alpha}) \in \mathfrak
g(\mathfrak o).$ The second condition is:
\[
\varepsilon^{k + \text{val}(\alpha(\ell))} Y_{\alpha} \in \mathfrak
g(\mathfrak o)
\]
or equivalently,  $\text{val}(\alpha(\ell)) + k \ge 0.$   But
(\ref{eqn-complement}) is a $\widetilde{T}(\mathbb C)$-invariant complement to
$\mathfrak s_{\ell}$ in $\mathfrak N(I)$ so the exponential mapping takes it
isomorphically to the Bruhat cell $I \ell x_0.$ \qed

\begin{lem}  The fixed point set of the turning torus on $X$ is:
\[ X^{\mathbb C^{\times}} = G(\mathbb C) \Lambda x_0.\]
The fixed point set of the torus $T(\mathbb C)$ on $X$ is:
\[ X^{T(\mathbb C)} = \Lambda x_0.\]
\end{lem}

\subsection{Proof}\label{subsec-fixedproof}
It is easy to see that $G(\mathbb C) \Lambda x_0 \subset X^{\mathbb
C^{\times}}.$  To verify the reverse inclusion it suffices to show that
\[ X^{\mathbb C^{\times}} \cap C_{\ell} \subset G(\mathbb C) \Lambda x_0 \]
for each Bruhat cell $C_{\ell}.$  But the fixed point set of the turning torus
on the vector space $D_{\ell}$ consists of those factors in the sum
(\ref{eqn-complement}) with $k=0.$ So their exponentials are contained in
$G(\mathbb C)\cap I.$  Similarly, the action of $T(\mathbb C)$ on $D_{\ell}$
has a single fixed point at the origin so $T(\mathbb C)$ acts on $C_{\ell}$
with a single fixed point $\ell x_0.$  \qed

\subsection{One dimensional orbits}
Let $X_1 \subset X$ denote the union of the 0 and 1-dimensional orbits under
the action of $T(\mathbb C)$ on $X.$  If $\alpha \in \Phi^+$ let $U_{\alpha}
\subset G$ be the unique connected $T$-stable 1-dimensional unipotent subgroup
whose Lie algebra contains the root space $\mathfrak g_{\alpha}.$  Then each
$\alpha \in \Phi^+$ determines a unique reductive connected algebraic subgroup
$H_{\alpha} \subset G$ of semi-simple rank 1 which is generated by $T,
U_{\alpha},$ and $U_{-\alpha}.$  Let $X^{\alpha} = H_{\alpha}(F)/H_{\alpha}
(\mathfrak o)$ be the loop Grassmannian for this subgroup.  The inclusion of
$H_{\alpha}$ into $G$ induces an injection $X^{\alpha} \hookrightarrow X.$

\begin{lem}\label{lem-1dimensional}
The union of the 0 and 1-dimensional orbits of $T(\mathbb C)$ on $X$ is
\[ X_1 = \bigcup_{\alpha \in \Phi^+} X^{\alpha}.\]
If $\alpha,\beta \in \Phi^+$ and $\alpha \ne \beta$ then $X^{\alpha} \cap
X^{\beta} = \Lambda.$
\end{lem}

\subsection{Proof} If $\dim(T)=1$ then the result is obvious since both sides
coincide with $X$, so we may assume that $\dim(T) \ge 2.$  The adjoint action
of the torus $T(\mathbb C)$ on $H_{\alpha}$ factors through an action of the
1-dimensional quotient $T(\mathbb C)/\ker(\alpha).$  Hence $X^{\alpha} \subset
X_1.$   To show that $X_1 \subset \bigcup_{\alpha\in \Phi^+} X^{\alpha}$ it
suffices to show that
\begin{equation}\label{eqn-X1piece} X_1 \cap C_{\ell} \subset
\bigcup_{\alpha\in
\Phi^+}X^{\alpha} \cap C_{\ell}\end{equation}
for each Bruhat cell $C_{\ell}.$  By Lemma \ref{lem-Bruhatcell} this amounts
to determining the 1-dimensional orbits of the action of $T(\mathbb C)$ on the
vector space $D_{\ell}$ of (\ref{eqn-complement}).  It follows from
(\ref{eqn-rootaction}) that if $\dim(T) \ge 2$ then these 1-dimensional orbits
are precisely the coordinate axes $\mathbb C \varepsilon^k Y_{\alpha}.$
Therefore we may classify the 1-dimensional orbits of $T(\mathbb C)$ on
$D_{\ell}$ by the positive roots, with those orbits corresponding to a fixed
$\alpha \in \Phi^+$ being contained in the subspace
\[D_{\ell,\alpha}= \bigoplus_k \mathbb C \varepsilon^k Y_{\alpha} =
\bigoplus_k \epsilon^k
\text{Lie}(U_{\alpha})(\mathbb C). \]
Here, the sum is taken over those $k$ such that $\alpha(a) + k > 0$ (for all
$a \in C_0$) and $\text{val}(\alpha(\ell)) + k < 0.$  Let $C_{\ell,\alpha}
\subset C_{\ell}$ be the corresponding subset of the Bruhat cell.  Then
$C_{\ell, \alpha} =  X^{\alpha}\cap C_{\ell},$ which implies
(\ref{eqn-X1piece}).  Finally, if $\alpha \ne \beta$ then $D_{\ell,\alpha}
\cap D_{\ell,\beta} = \left\{ 0 \right\}$ which proves that $X^{\alpha} \cap
X^{\beta} = \Lambda.$  \qed

\section{Affine Springer fibers for $\text{SL}(2)$}  \label{sec-SL2}
\subsection{}
In this section we fix $G(\mathbb C) = \text{SL}(2, \mathbb C)$ and set
$X^{\SL} = G(F)/K$ with basepoint  $x_0=K=G(\mathfrak o).$  Let $T$ be the
torus of diagonal matrices and let $\alpha$ be the positive root, $\alpha
\left( \begin{smallmatrix} a & 0 \\ 0 & a^{-1} \end{smallmatrix}\right) =
a^2.$  Let $\alpha^{\vee}: F^{\times} \to T(F)$ be the corresponding co-root,
\[ \alpha^{\vee}(b) = \left( \begin{matrix} b & 0 \\ 0 & b^{-1} \end{matrix}
\right).\]
Then $\Lambda^{\SL} = \langle \alpha^{\vee}(\epsilon) \rangle$ is the
cyclic group spanned by $\alpha^{\vee}(\epsilon),$ whose elements we denote by
$\ell_m =
\alpha^{\vee}(\epsilon^m).$  For any integer $n \le -1$ set
\begin{equation}\label{eqn-basepoint}
  x_n = \left(\begin{matrix} 1 & \epsilon^{n} \\ 0 & 1 \end{matrix} \right)K
\in X^{\SL}.\end{equation}
Let $\alpha':\mathfrak t(F) \to F$ be the differential of $\alpha:T(F) \to F.$
The following fact was observed by D. Nadler.

\begin{lem}\label{lem-Nadler}
The affine Grassmannian $X^{\SL}$ is the disjoint union of countably many
$T(F)$ orbits,
\[ X^{\SL} = \coprod_{n \le 0} T(F) \cdot x_n.\]
The orbit of the point $x_n$ has \parens{complex} dimension $\abs{n}.$   If
$\gamma \in \mathfrak t(\mathfrak o)$ and $v = \val(\alpha'(\gamma))$ then the
affine Springer fiber $X_{\gamma}^{\SL}$ is the union of the $T(F)$ orbits
\begin{equation}\label{eqn-unionoforbits}
  X_{\gamma}^{\SL} = \bigcup_{n=-v}^0 T(F)\cdot x_n  \end{equation}
and it is preserved by the turning torus.
\end{lem}
\noindent
It follows that we may unambiguously denote the affine Springer fiber
$X_{\gamma}^{\SL}$ in \parens{\ref{eqn-unionoforbits}} by $X^{\SL}_{\le v}.$

\subsection{Proof}
First observe that every point $x \in X^{\SL}$ has an expression of the form
$gK$ with
\[ g = \left( \begin{matrix} \epsilon^m & b_0 \epsilon^n \\ 0 & \epsilon^{-m}
\end{matrix}
\right) \]
where either (1) $b_0 = 0$ or (2) $b_0 \in \mathfrak o^{\times}$ and $n-m <
0.$  In fact, the standard Borel subgroup acts transitively on $X,$ so  $x =
\left( \begin{smallmatrix} a & b' \\ 0 & a^{-1} \end{smallmatrix} \right)x_0$
for some $a,b' \in F.$  Set $a = a_0 \epsilon^m$ with $a_0 \in \mathfrak
o^{\times}$ and right multiply by $\alpha^{\vee}(a_0^{-1}) \in K$ to obtain
$x=gx_0$ with $g = \left( \begin{smallmatrix} \epsilon^m & b \\ 0 &
\epsilon^{-m} \end{smallmatrix} \right).$  If $b =0$ or if $\val (b) < m$ then
we are done.  Otherwise set $b = b_0 \epsilon^n$ with $n \ge m$ and $b_0 \in
\mathfrak o^{\times}.$  Right multiply by $k = \left( \begin{smallmatrix} 1 &
-b_0 \epsilon^{n-m} \\ 0 & 1 \end{smallmatrix} \right) \in K$ to see that $x =
gkK\in \Lambda.$

Now let $x = gK$ with $g = \left( \begin{smallmatrix} \epsilon^m & b_0
\epsilon^n \\ 0 & \epsilon^{-m} \end{smallmatrix} \right).$  If $b_0 = 0$ then
$x = \alpha^{\vee}(\epsilon^m)x_0 \in T(F)\cdot x_0.$  If $b_0 \in \mathfrak
o^{\times}$ and $m>n$ let $a \in \mathfrak o^{\times}$ be a square root of
$b_0.$  Set $ t = \alpha^{\vee}(a \epsilon^m) \in T(F)$ and $k = \alpha^{\vee}
(a^{-1}) \in K.$  Then $tx_{n-m}k = g$ which proves the first statement.  The
dimension statement is straightforward. If $\gamma \in \mathfrak t(\mathfrak
o)$ then $X_{\gamma}^{\SL}$ is preserved by $T(F)$ so it is a union of $T(F)$
orbits.  To verify (\ref{eqn-unionoforbits}), suppose $\gamma = \left(
\begin{smallmatrix} a & 0 \\ 0 & -a \end{smallmatrix} \right) \in \mathfrak
t(\mathfrak o).$  Then $\text{Ad}(x_n)\gamma = \left( \begin{smallmatrix} a &
-2a\epsilon^n\\ 0 & -a \end{smallmatrix} \right)$ which is in $\mathfrak
g(\mathfrak o)$ iff $\val(a) +n \ge 0.$  \qed

Each 1-dimensional orbit $\mathcal O$  of $\widetilde{T}(\mathbb C)$
necessarily contains exactly two fixed points $\ell_s,\ell_t\in\Lambda$ in its
closure, in which case we say that $\mathcal O$ ``connects'' the points
$\ell_s$ and $\ell_t.$

\begin{lem}\label{lem-SL2orbits} For any two lattice points $\ell_s, \ell_t
\in \Lambda$ there exists a unique 1-dimensional orbit $\mathcal O_{st}$ of
$\widetilde{T}(\mathbb C)$ in $X^{\SL}$ which connects them, and this accounts
for all the 1-dimensional orbits of $\widetilde{T}(\mathbb C)$ in $X^{\SL}.$
The stabilizer of any point in $\mathcal O_{st}$ is the kernel of the affine
root $(\alpha, s+t).$  The orbit $\mathcal O_{st}$ is contained in $X_{\gamma}
^{\SL}$ iff $\val(\alpha'(\gamma)) \ge \abs{s-t}.$
\end{lem}

\subsection{Proof}  By Lemma \ref{lem-Nadler} for any $T(\mathfrak o)$ orbit
$\mathcal T$ there exists a unique $s \in \mathbb Z$ and a unique $n \le 0$
such that $\ell_s x_n \in \mathcal T.$  Set $t = s+n.$  We will show that
\begin{enumerate}
\item the $\widetilde{T}(\mathbb C)$ orbit $\mathcal O$ of $\ell_s x_n$ is
1-dimensional,
\item this is the only 1-dimensional orbit of $\widetilde{T}(\mathbb C)$ in
$\mathcal T,$
\item $\mathcal O$ connects $\ell_s$ and $\ell_t=\ell_s + n \alpha^{\vee}
(\epsilon),$ and
\item $\widetilde{T}(\mathbb C)$ acts on $\mathcal O$ through the affine root
$(\alpha, 2s+n).$
\end{enumerate}

First compute the action of $\widetilde{T}(\mathbb C)$ on the $T(\mathfrak o)$
orbit of the point $\ell x_n.$ Let $\alpha^{\vee}(b) \in T(\mathfrak o)$ with
$ b = \sum_{i \ge o} b_i \epsilon^i \in \mathfrak o^{\times}.$  Let $a \in
\mathbb C^{\times}.$  Since
\[ \ell_s x_n = \left( \begin{matrix} \epsilon^s & \epsilon^{s+n} \\
0 & \epsilon^{-s} \end{matrix} \right) \]
the action of  $(\alpha^{\vee}(a),\lambda) \in \widetilde{T}(\mathbb C)$ on
the point $\alpha^{\vee}(b) \ell_s x_n$ is the point
\[ y = \left( \begin{matrix}
a(\lambda \cdot b) \lambda^s \epsilon^s & a (\lambda\cdot b) \lambda^{s+n}
\epsilon^{s+n} \\
0 & a^{-1} (\lambda \cdot b^{-1}) \lambda^{-s} \epsilon^{-s} \end{matrix}
\right) K. \]
Right multiply by $\alpha^{\vee}( a^{-1} (\lambda \cdot b)^{-1} \lambda^{-s})
\in K$ to find that
\begin{equation}\label{eqn-1daction}
  y = \left( \begin{matrix} \epsilon^s & a^2 (\lambda \cdot b)^2 \lambda^{2s+n}
\epsilon^{s+n} \\
0 & \epsilon^{-s} \end{matrix} \right) K .\end{equation}
As $a, \lambda \in \mathbb C^{\times}$ vary, this is a 2-dimensional orbit
unless $b = b_0 \in \mathbb C^{\times}$, in which case $\lambda \cdot b = b,$
hence $\lim_{a\to 0} y = \ell_s.$ To find $\lim_{a \to \infty} y$, choose
$\lambda = b = 1$ and right multiply by
\[ \left( \begin{matrix} a^{-3} & a^{-1} \epsilon^{-n} - a^2 \\
a^{-2} & \epsilon^{-n} \end{matrix} \right) \in K\]
to find
\[ y = \left( \begin{matrix}
a^{-3} \epsilon^s + \epsilon^{s+n} & a^{-1} \epsilon^{s-n} \\
a^{-2} \epsilon^{-s} & \epsilon^{-s-n} \end{matrix} \right)K. \]
So $y \to \ell_s + n\alpha^{\vee}(\epsilon)$ as $a \to \infty.$ It follows
from (\ref{eqn-1daction}) that $\widetilde{T}(\mathbb C)$ acts on this
1-dimensional orbit through the character $(\alpha, 2s+n)=(\alpha, s+t)$, cf.
(\ref{eqn-rootaction}). This verifies (1)--(4).  Finally, by
(\ref{eqn-unionoforbits}), the orbit $\mathcal O \subset X_{\gamma}$ iff
$\abs{n} \le \val \alpha'(\gamma).$  \qed

\section{Equivariant homology of $\text{SL}(2)$ Springer fibers}
\label{sec-EqHoSL2}
\subsection{}
As in \S \ref{sec-SL2} let $G = \text{SL}(2),$ let $X = X^{\text{SL}(2)} =
G(K)/G(\mathfrak o)$, let $\Lambda = \Lambda^{\text{SL}(2)}.$  (We use
multiplicative notation for the group law in $\Lambda.$)  Let $T\subset G$ be
the diagonal matrices, with simple root and co-root $\alpha$ and
$\alpha^{\vee}.$  Then $\alpha^{\vee}(\epsilon)\in \Lambda$ is a generator
(which, by abuse of notation, we denote by $\alpha^{\vee}$) so it determines
an isomorphism $\mathbb Z \cong \Lambda,$ written $s \mapsto \ell_s.$ So
$\alpha^{\vee} = \ell_1$ and  $\alpha^{\vee}\ell_s = \ell_{s+1}.$ Fix $\gamma
\in \mathfrak t(\mathfrak o)$ with corresponding Springer fiber $X_{\gamma}.$

To simplify the notation, for the remainder of \S \ref{sec-EqHoSL2} we will
write $\widetilde{T}$ for $\widetilde{T}(\mathbb C)$, $\mathfrak t$ for
$\mathfrak t(\mathbb C)$, and so on.  Let $\mathcal D = \mathcal D(\mathfrak
t)$  be the algebra of differential operators on $\mathfrak t^*,$ which is
identified with $H^*_T(\text{pt})$ by the Chern-Weil homomorphism.  Let
$\mathbf S(\mathfrak t)$ be the symmetric algebra of complex valued polynomial
functions on $\mathfrak t^*.$  The character $\alpha$ corresponds to a
differential operator $\partial_{\alpha} \in \mathcal D(\mathfrak t).$  Denote
by $\bS(\mathfrak t)\left\{ \partial_{\alpha}^d \right\}$ the kernel of
$\partial_{\alpha} ^d.$  The Chern-Weil homomorphism determines an isomorphism
of $\mathcal D$ modules
\[ H^{T}_*(\Lambda) \cong \mathbb C[\Lambda] \otimes_{\mathbb C} \mathbb
\bS(\mathfrak t). \]

\begin{prop}\label{prop-SL2answer}  Let $v = \val(\alpha'(\gamma)).$  The
inclusion $\Lambda \subset X$ induces  exact sequences,
\begin{equation}\label{eqn-maindiag}\begin{CD}
H_*^{\widetilde{T}}(X_{\gamma},\Lambda) @>{\tilde\beta}>> H_*^{\widetilde{T}}
(\Lambda) @>>>
H_*^{\widetilde{T}}(X_{\gamma}) @>>> 0 \\
\bigcup && \bigcup && \bigcup \\
H_*^{T}(X_{\gamma},\Lambda) @>{\beta}>> H_*^{T}(\Lambda) @>>>
H_*^{T}(X_{\gamma}) @>>> 0
\end{CD}\end{equation}
and the image of $\beta$ is the $\mathcal D$-submodule
\begin{equation}\label{eqn-imageofbeta}
\sum_{d =1}^v (1-\alpha^{\vee})^d \mathbb C[\Lambda] \otimes \bS(\mathfrak t)
\left\{\partial^d_{\alpha}\right\} \subset \mathbb C[\Lambda] \otimes
\bS(\mathfrak t).
\end{equation}
\end{prop}

\subsection{Proof}
The Springer fiber $X_{\gamma}$ has a paving by affines (\cite{GKM2}):
it is an increasing union
of complex projective algebraic varieties $\phi \subset
(X_{\gamma})_0 \subset (X_{\gamma})_1
\subset \ldots$ such that each $(X_{\gamma})_j - (X_{\gamma})_{j-1}$
is isomorphic to a single
affine space. In fact the intersection of $X_{\gamma}$ with a single
Bruhat cell is either empty,
or else it is a single affine cell in this paving.  It follows that
the homology of
$(X_{\gamma})_j$ vanishes in odd degrees.  The extended torus
$\widetilde{T} = T \times \mathbb
C^{\times}$ (\S \ref{subsec-turning}) preserves each $(X_{\gamma})_j$ and it
acts with finitely many fixed points and finitely many 1-dimensional orbits.
So the results of \S \ref{subsec-equivariant1} and \S
\ref{subsec-equivariant2} may be applied to this action. Since homology
commutes with direct limits, we obtain the exact sequences
(\ref{eqn-maindiag}) for both the $\widetilde{T}$-equivariant homology and the
$T$-equivariant homology.

The co-root $\alpha^{\vee}$ determines  canonical isomorphisms
\begin{equation}\label{eqn-Cxt}
  \bS(\widetilde{\mathfrak t}) \cong \mathbb C[x,t] \ \text{ and } \ \mathcal
D(\widetilde{\mathfrak t}) \cong \mathbb
C[\partial_x,\partial_t].\end{equation}
According to Lemma \ref{lem-SL2orbits}, for each unordered pair of distinct
integers $a,b \in \mathbb Z$ such that $\abs{b-a} \le v=
\val(\alpha'(\gamma)),$ there is a unique 1-dimensional orbit $E_{ab}\subset
X_{\gamma}$ which connects $\ell_a$ and $\ell_b$, on which the torus
$\widetilde{T}$ acts through the character $\Phi_{ab}$ whose differential
$\phi_{ab}: \widetilde{\mathfrak t} \to \mathbb C$ corresponds to the
differential operator
\[ \partial_{ab} = \partial_{\alpha} + (a+b)\partial_t = 2 \partial_x  + (a+b)
\partial_t.\]
Let $\mathfrak m_{ab} = \ker(\phi_{ab})$ denote the Lie algebra of the
stabilizer of any point in this orbit, and let $j_{ab}^*:\bS(\mathfrak m_{ab})
\to \mathcal S(\widetilde{\mathfrak t})$ be the resulting inclusion.  Then by
(\ref{eqn-annihilator}) the image of $j_{ab}^*$ consists of all polynomial
functions on $\widetilde{\mathfrak t}^*$ which are annihilated by the
differential operator $\partial_{ab}.$ Hence
\[
j_{ab}^*\bS(\mathfrak m_{ab}) = \left\{ g((a+b)x - 2t):\ g \in \mathbb C[z]
\right\}
\]
consists of polynomial functions of $z=(a+b)x-2t.$  It follows that the image
of $\tilde\beta$ is the sum $\sum_{a,b} M_{ab}$ of submodules $M_{ab}$ spanned
by elements
\[  (\ell_b-\ell_a)\otimes g_{ab}((a+b)x-2t)\]
with $\abs{b-a} \le v.$  This is the module $P_v$ of \S \ref{subsec-Pv} (with
$2t$ replaced by $t$).

The $T$-equivariant homology of $X_{\gamma}$ may be recovered
(\ref{eqn-annihilator}) from the $\widetilde{T}$-equivariant homology of
$X_{\gamma}$ as the kernel of the operator $\partial_t.$ So the image of
$\beta$ is
\[
\text{Im}(\beta) = P_v \cap \ker(\partial_t).
\]
Then Proposition \ref{prop-Fmdspan} identifies $P_v\cap \ker(\partial_t)$ with
the submodule (\ref{eqn-imageofbeta}). \qed

\section{Groups of semisimple rank one}\label{sec-ssrankone}
For lack of an adequate reference we include a proof of the following
well-known fact.
\begin{lem} \label{lem-classification}
Let $H$ be a connected reductive complex linear algebraic group of rank $n$
and of semisimple rank 1.  Then $H$ is isomorphic to either \begin{enumerate}
\item $( \mathbb C^{\times})^{n-1} \times \text{\rm SL}(2, \mathbb C)$
\item $ (\mathbb C^{\times})^{n-1} \times \text{\rm PGL}(2, \mathbb C)$
\item $ (\mathbb C^{\times})^{n-2} \times \text{\rm GL}(2, \mathbb C).$
\end{enumerate}\end{lem}
\subsection{Proof}
Let $(L,\left\{ \alpha \right\}, L^{\vee}, \left\{ \alpha^{\vee} \right\})$ be
the (based) root datum for $H$ where $\langle \cdot, \cdot \rangle: L^{\vee}
\times L \to \mathbb Z$ is a dual pairing of rank $n$ lattices, $\alpha \in
L$, $\alpha^{\vee} \in L^{\vee}$, and $\langle \alpha^{\vee}, \alpha \rangle =
2.$  If
\begin{equation}\label{eqn-function} \ \langle \alpha^{\vee}, \cdot \rangle :
L \to \mathbb Z\end{equation}
is even-valued, then $\frac{1}{2}\alpha^{\vee} \in L^{\vee}$ so
\[ L = \left( \alpha^{\vee} \right) ^{\perp} \oplus \langle \alpha \rangle
\text{\ while\ } L^{\vee} = \alpha^{\perp} \oplus \langle \textstyle{\frac{1}
{2}} \alpha^{\vee} \rangle\]
which corresponds to case (1).  (Here, $\langle \alpha \rangle$ denotes the
cyclic group generated by $\alpha \in L,$ and $(\alpha^{\vee})^{\perp}$
denotes the kernel of (\ref{eqn-function}).)  Similarly if $\langle \cdot,
\alpha \rangle : L^{\vee} \to \mathbb Z$
is even-valued then $\frac{1}{2}\alpha \in L$ so
\[L = \left( \alpha^{\vee} \right)^{\perp} \oplus \langle \textstyle{\frac{1}
{2}} \alpha \rangle \text{\ while\ } L^{\vee} = \alpha^{\perp} \oplus \langle
\alpha^{\vee} \rangle\]
which corresponds to case (2).

Now suppose that neither of these homomorphisms is even-valued.  Let $\left\{
e_1,e_2, \ldots, e_{n-1}\right\}$ be a basis of $L_0 = \left( \alpha^{\vee}
\right)^{\perp}.$  Let $L_1 = \left\{ x\in L :\ \langle \alpha^{\vee},x
\rangle  = 1 \right\}.$  Then $L_1-\frac{1}{2} \alpha$ is a translate of $L_0$
so there exists $a_i \in \left\{ 0, \frac{1}{2} \right\}$ such that $L_1 -
\frac{1}{2} \alpha = L_0 + \sum_{i=1}^{n-1}a_ie_i$ and not all the $a_i$ are
0.  By relabeling the basis we may assume that $a_1 = a_2 = \ldots = a_r =
\frac{1}{2}$ and the remaining coefficients are 0.  Let $v = e_1 + e_2 +
\ldots + e_r.$ Then $\left\{ v, e_2, \ldots, e_r\right\}$ is also a basis for
$L_0,$ and $\frac{1}{2} \alpha +\frac{1}{2} v \in L.$  It follows that
\[ L = Y \oplus \langle e_2,e_3, \ldots, e_{n-1} \rangle \]
where
\[ Y = \left\{ m\alpha + m'v :\ m,m' \in \mathbb Z \text{\ or\ } m,m' \in
\mathbb Z + \textstyle{\frac{1}{2}}  \right\}. \]
Then this corresponds to case (3), with $L^{\vee} = \langle
e_2,e_3,\ldots,e_{n-1} \rangle^{\perp} \oplus Y^{\perp}.$  \qed

\subsection{}\label{subsec-H1}
Let $H$ be a connected reductive complex linear algebraic group of semisimple
rank one.  Fix a Borel pair $B \subset T \subset H$ with resulting positive
root $\alpha.$  Let $X^H$ be the affine Grassmannian for $H$ and let
$\Lambda^H \subset X^H$ be the lattice of translations for $T.$ Consider the
abstract $\text{SL}(2,\mathbb C)$ which corresponds to the simple root
$\alpha.$  The canonical mapping $\text{SL}(2) \to H$  induces an inclusion
$\phi_{\alpha}: X^{\text{SL}(2)} \subset X^H$ with $\phi_{\alpha}
(\Lambda^{\SL}) = \langle \alpha^{\vee}(\epsilon) \rangle \subset \Lambda^H$
(the infinite cyclic subgroup generated by $\alpha^{\vee}(\epsilon)$). To
simplify the notation, identify the co-root $\alpha^{\vee}$ with its image
$\alpha^{\vee}(\epsilon)$ in $\Lambda^H.$

\begin{lem} \label{lem-SS1}
The following statements hold. \begin{enumerate}
\item The affine Grassmannian for $H$ is the disjoint union,
\[ X^H = \coprod_{\ell \in \Lambda^H/\langle \alpha^{\vee} \rangle } \ell
\phi_{\alpha}(X^{\SL}).\]
\item  The $T^H(\mathfrak o)$ orbits on $X^H$ coincide with the $T^{\SL}
(\mathfrak o)$ orbits on $X^H$ \parens{which are in turn translates, by
elements of $\Lambda^H$, of $T^{\SL}(\mathfrak o)$ orbits on $X^{\SL}$}.
\item  Let $\gamma \in \mathfrak t^H(\mathfrak o)$ be a regular element.  Then
the affine Springer fiber $X^H_{\gamma}$ is the disjoint union
\begin{equation}\label{eqn-Springerunion}
X^H_{\gamma} = \coprod_{\ell \in \Lambda^H/ \langle \alpha^{\vee} \rangle}
\ell \phi_{\alpha}(X^{\SL}_{\le v})
\end{equation}
where $ v = \text{\rm val}(\alpha'(\gamma)),$ cf. equation
\parens{{\rm\ref{eqn-unionoforbits}}}.
\end{enumerate}\end{lem}

\subsection{Proof}
For part (1) use Lemma \ref{lem-classification} to reduce to the case $H=
\text{GL} (2)$ or $\text{PGL}(2).$  If $H=\text{GL}(2)$ then the exact
sequences
\begin{equation*}\begin{CD}
I @>>> \text{SL}(2,F) @>>> \text{GL}(2,F) @>{\text{det}}>> F^{\times} @>>> 1
\\
I @>>> \text{SL}(2, \mathfrak o) @>>> \text{GL}(2, \mathfrak o) @>>> \mathfrak
o^{\times} @>>> 1 \end{CD} \end{equation*}
together with the isomorphism $\val: F^{\times} / \mathfrak o^{\times} \cong
\mathbb Z$ give rise to a diagram
\begin{equation*} \begin{CD}
&& X^{\SL} @>>> X^{\text{GL}(2)} @>{\val\det{}}>> \mathbb Z \\
&& @AAA @AAA @| \\
0 @>>> \Lambda^{\SL} @>>> \Lambda^{\text{GL}(2)} @>>> \mathbb Z @>>> 0
\end{CD} \end{equation*}
from which the result follows.  Now suppose $H = \text{PGL}(2).$  If $g \in
H(F)$ let $\tau(g) = \val\det{g} \mod 2.$  If $\tau(g)=0$ then $\det g$ has a
square root in $F$ and $g/\sqrt{\det g} \in \text{SL}(2,F).$ So the following
exact sequences
\begin{equation*} \begin{CD}
\left\{ \pm I \right\} @>>> \text{SL}(2,F) @>>> \text{PGL}(2,F) @>{\tau}>>
\mathbb Z/(2)\\
\left\{ \pm I \right\} @>>> \text{SL}(2,\mathfrak o) @>>> \text{PGL}
(2,\mathfrak o) @>>> 1
\end{CD} \end{equation*}
give rise to a diagram
\begin{equation*} \begin{CD}
&& X^{\SL} @>>> X^{\text{PGL}(2)} @>{\tau}>> \mathbb Z/(2)\\
&& @AAA @AAA @| \\
0 @>>> \Lambda^{\SL} @>>> \Lambda^{\text{PGL}(2)} @>>> \mathbb Z/(2) @>>> 0
\end{CD} \end{equation*}
from which the result follows.

For part (2) we may assume $H = \text{GL}(2)$ or $\text{PGL}(2).$  Let $gK \in
X^H$ and let $t =\left(\begin{smallmatrix} a & 0 \\ 0 & d \end{smallmatrix}
\right) \in T^H(\mathfrak o).$  Then $ad^{-1}$ has a square root $\lambda \in
\mathfrak o.$  Set $t' = \alpha^{\vee}(\lambda) \in T^{\SL}(\mathfrak o).$  By
matrix multiplication, $tgt^{-1} = t' g (t')^{-1}$ hence $tgK^H = t' gK^H.$

Now consider part (3).  It is easy to verify for $H=\text{GL}(2),$  for $H =
\text{PGL}(2)$ and for $n \le -1$ that $x_n \in X^H_{\gamma}$ iff $\abs{n} \le
v = \val(\alpha'(\gamma)),$ where $x_n = \left(\begin{smallmatrix} 1 &
\epsilon^n \\ 0 & 1 \end{smallmatrix} \right)K.$ Since $X_{\gamma}^H$ is a
union of $T^{\SL}(\mathfrak o)$ orbits, it follows that
\begin{equation}\label{eqn-intersect}
X_{\gamma}^H \cap \phi_{\alpha}(X^{\SL}) = \phi_{\alpha}(X^{\SL}_{\le v}).
\end{equation}
By Lemma \ref{lem-classification}, equation (\ref{eqn-intersect}) holds for
any $H$ of semi-simple rank one. It follows that the union
(\ref{eqn-Springerunion}) is contained in $X_{\gamma}^H.$ Conversely, if
$\mathcal O$ is a $T^H(\mathfrak o)$ orbit in $X_{\gamma}^H$ then there exists
$\ell \in
\Lambda^H$ so that $\ell \phi_{\alpha}(x_n) \in \mathcal O,$ which holds iff $
\phi_{\alpha}(x_n) \in X_{\gamma}^H $ or $x_n \in X^{\SL}_{\le v}.$  So
\begin{equation*}
X_{\gamma}^H \subset \bigcup_{\ell \in \Lambda^H} \ell \phi_{\alpha}
\left(X^{\SL}_{\le v} \right)
\end{equation*}
The normalizer of $\phi_{\alpha}(X^{\SL}_{\le v})$ in $\Lambda^H$ is the
sub-lattice $\langle \alpha^{\vee} \rangle = \phi_{\alpha}(\Lambda^{\SL})$, so
the union (\ref{eqn-Springerunion}) is disjoint.  \qed

\section{Equivariant homology of affine Springer fibers}
\label{sec-equivhomology}
\subsection{}
As in \S \ref{sec-affinespringer} we consider a connected reductive linear
algebraic group $G$ defined over $\mathbb C$ and a Borel pair $T\subset B
\subset G$ with resulting system of positive roots $\Phi^+.$  Fix $\gamma \in
\mathfrak t(\mathfrak o)$ and let $X_{\gamma}$ be the resulting affine
Springer fiber.  It contains the lattice $\Lambda = \Lambda^G.$  We use
multiplicative notation for the group operation in $\Lambda.$  For notational
simplicity in this section we denote by $\mathfrak t, \mathfrak t^*$, etc.,
the complex vector spaces $\mathfrak t(\mathbb C)$, $\mathfrak t^*(\mathbb C)
$,
etc.  Let $\alpha \in \Phi^+$ and let $\alpha^{\vee}:F^{\times} \to T(F)$ be
the corresponding co-root; it determines a lattice element which we also
denote by $\alpha^{\vee} \in \Lambda.$ Moreover, $\alpha^{\vee} \in \mathfrak
t = \mathfrak t^{**}$ corresponds to a (degree 1) monomial
\[ x_{\alpha^{\vee}} \in \bS(\mathfrak t) \]
while $\alpha \in \mathfrak t^*$ corresponds to a (degree 1) differential
operator
\[ \partial_{\alpha} \in \mathcal D(\mathfrak t) \]
such that $\partial_{\alpha}(x_{\alpha^{\vee}}) = 2.$  Let $\bS(\mathfrak t)
\left\{\partial_{\alpha}^d\right\}$ be the submodule of polynomial functions
which are annihilated by the differential operator $\partial_{\alpha}^d.$
Define the following submodule of $H^{T(\mathbb C)}_*(\Lambda) = \mathbb
C[\Lambda] \otimes \bS(\mathfrak t):$
\[ L_{\alpha, \gamma} = \sum_{d=1}^{\text{val} (\alpha'(\gamma))}
\left(1-\alpha^{\vee}\right)^d    \mathbb C[\Lambda] \otimes
\bS(\mathfrak t) \left\{\partial_{\alpha}^d\right\}.\]

\begin{thm}\label{thm-cohomology}  Suppose the homology $H_*(X_{\gamma};
\mathbb C)$ is pure. Then the inclusion $\Lambda \subset X_{\gamma}$ induces
an exact sequence
\begin{equation*} \begin{CD}
0 @>>> \underset{\alpha \in \Phi^+}{\sum} L_{\alpha,\gamma} @>>>
\mathbb C[\Lambda] \otimes \bS(\mathfrak t) @>>> H_*^{T(\mathbb C)}
(X_{\gamma}) @>>> 0.
\end{CD}\end{equation*}
\end{thm}

\subsection{Proof}\label{subsec-proof-cohomology}
Since the homology of $X_{\gamma}$ is equivariantly formal, the exact sequence
(\ref{eqn-localization}) becomes
\begin{equation} \label{eqn-SS1}
\begin{CD}
H^{T(\mathbb C)}_*\left((X_{\gamma}\right)_1, \Lambda) @>{\beta}>>
H^{T(\mathbb C)}_*(\Lambda) @>>> H^{T(\mathbb C)}_*(X_{\gamma}) @>>> 0
\end{CD} \end{equation}
where $(X_{\gamma})_1$ denotes the closure of the set of 1-dimensional
$T(\mathbb C)$ orbits in $X_{\gamma}.$  By Lemma \ref{lem-1dimensional},
$X^{\alpha}_{\gamma} \cap X^{\beta}_{\gamma} = \Lambda$ if $\alpha \ne \beta,$
and
\[ \left( X_{\gamma}\right)_1 = \bigcup_{\alpha \in \Phi^+}
X^{\alpha}_{\gamma} \]
where $X^{\alpha}_{\gamma}$ is the affine Springer fiber corresponding to
$\gamma $ in the loop Grassmannian $X^{\alpha}$ for the group $H^{\alpha}$ of
semi-simple rank one which is determined by $\alpha.$  So the image of $\beta$
is the sum over positive roots of the image of the corresponding mapping in
the semisimple rank one case. Thus, it suffices to consider the case in which
the group $G$ has semisimple rank one, which we now assume.

Let $A = \ker(\alpha)^0 \subset T$ and let $T_1 \subset T$ be the
1-dimensional sub-torus corresponding to the co-root $\alpha^{\vee}.$  The
canonical decomposition $\mathfrak t = \mathfrak t_1 \oplus \mathfrak a$
determines an isomorphism $H^{T(\mathbb C)}_*(\text{pt})\cong H^{T_1(\mathbb
C)}(\text{pt}) \otimes H^{A(\mathbb C)}_*(\text{pt})$, that is,
$\bS(\mathfrak t)  \cong \bS(\mathfrak t_1) \otimes \bS(\mathfrak a) $ and
\begin{equation*}
\bS(\mathfrak t)\left\{\partial^d_{\alpha}\right\} \cong
\bS(\mathfrak t_1)\left\{\partial^d_{\alpha}\right\} \otimes \bS(\mathfrak a).
\end{equation*}

So the exact sequence (\ref{eqn-SS1}) for $T(\mathbb C)$-equivariant homology
is obtained from the same exact sequence for $T_1(\mathbb C)$-equivariant
homology
\begin{equation}\label{eqn-beta1}
\begin{CD}
H_*^{T_1(\mathbb C)}(X_{\gamma},\Lambda) @>{\beta_1}>> H_*^{T_1(\mathbb C)}
(\Lambda) @>>> H_*^{T_1(\mathbb C)}(X_{\gamma}) @>>> 0
\end{CD}\end{equation}
by tensoring with $\bS(\mathfrak a).$  By Lemma \ref{lem-SS1} the homomorphism
$\beta_1$ decomposes as a direct sum of homomorphisms
\begin{equation*}\begin{CD}
\underset{\ell\in \Lambda^H/ \langle \alpha^{\vee} \rangle}{\bigoplus}
\ell_* \phi_{\alpha *} H_*^{T_1(\mathbb C)}(X^{\SL}_{\le v}, \Lambda^{\SL})
@>>> \underset{\ell\in \Lambda^H/ \langle \alpha^{\vee} \rangle}{\bigoplus}
\ell_* \phi_{\alpha *} H_*^{T_1(\mathbb C)}(\Lambda^{\SL})
\end{CD} \end{equation*}
(where $\phi_{\alpha*}$ and $\ell_*$denote the homomorphisms on homology which
are induced by $\phi_{\alpha}$ and by translation by $\ell \in \Lambda^H,$
respectively).  So by Proposition \ref{prop-SL2answer} the image of $\beta_1$
is the sum
\begin{equation*}
\sum_{\ell\in \Lambda^H/ \Lambda^{\SL}} \ell_* \phi_{\alpha *}\sum_{d=1}
^v(\alpha^{\vee}-1)^d \mathbb C[\Lambda^{\SL}] \otimes \bS(\mathfrak t_1)
\left\{\partial^d_{\alpha}\right\}
= \sum_{d=1}^v(\alpha^{\vee}-1)^d \mathbb C[\Lambda] \otimes \bS(\mathfrak
t_1) \left\{\partial_{\alpha}^d\right\}.
\end{equation*}
Since the image of $\beta$  is obtained by tensoring with $\bS(\mathfrak a),$
we obtain
\[ \text{Im}(\beta) = \sum_{d=1}^v(\alpha^{\vee}-1)^d \mathbb C[\Lambda]
\otimes \bS(\mathfrak t) \left\{\partial_{\alpha}^d\right\}.  \qed\]

In fact, in the semisimple rank one case, as a module over $\bS(\mathfrak a)$
the vector space of relations $L_{\alpha, \gamma}$ has a basis consisting of
the collection of elements
\[ f_{\ell, d, \alpha} = (1-\alpha^{\vee})^d \ell \otimes x_{\alpha}^{d-1} \in
\mathbb C[\Lambda] \otimes \bS(\mathfrak t)\]
with $\ell \in \Lambda$ and $1 \le d \le \text{val}(\alpha'(\gamma)).$

\subsection{}\label{subsec-Weylaction}
Let $\mathbf{Aut}$ denote the automorphism group of the based root
datum for $G.$
Let $W = W(G, T)$ be the Weyl group and let $\widetilde{W} = \Lambda \rtimes
W$ be the extended affine Weyl group.  Then the group
\begin{equation}\label{eqn-automorphismgroup}
  \widetilde{W} \rtimes \mathbf{Aut} \end{equation}
acts on $T$ and on $\Lambda.$  It acts through $W \rtimes \mathbf{Aut}$ on
$T(\mathbb C)$ and on the root system $\Phi(G,T).$  Consider the  diagonal
action (from the left) of the group (\ref{eqn-automorphismgroup}) on the
equivariant homology
\[ H^{T(\mathbb C)}_*(\Lambda) \cong \mathbb C[\Lambda] \otimes_{\mathbb C}
\mathbf{S} (\mathfrak t).\]
It preserves the kernel
\[ L_{\infty} = \sum_{\alpha \in \Phi^+} \sum_{d=1}^{\infty} (1-\alpha^{\vee})
^d
\mathbb C[\Lambda] \otimes \bS(\mathfrak t)\left\{\partial_{\alpha}^d\right\}
\]
of the surjection $H^{T(\mathbb C)}_*(\Lambda) \to H^{T(\mathbb C)}_*(X)$ so
it passes to an action on the equivariant homology of $X,$ which we will refer
to as the ``left'' action. If $\gamma \in \mathfrak t(\mathfrak o)$ then the
subgroup
\begin{equation}\label{eqn-Wgamma}
\left( \widetilde{W} \rtimes \mathbf{Aut} \right)_{\gamma} = \left\{ \tau \in
\widetilde{W} \rtimes \mathbf{Aut}: \ \val(\tau\alpha(\gamma)) =
\val(\alpha(\gamma)) \text{
for all }\alpha \in \Phi^+ \right\}\end{equation}
preserves the module of relations
\[ \sum_{\alpha \in \Phi^+} L_{\alpha, \gamma} \subset H_*^{T(\mathbb C)}
(\Lambda).\]
Consequently, if the homology $H_*(X_{\gamma})$ is pure, the group
$(\widetilde W \rtimes \mathbf{Aut})_{\gamma}$ acts (from the left) on the
equivariant homology $H^{T(\mathbb C)}(X_{\gamma})$ by $\mathcal D(\mathfrak
t)$ module homomorphisms.  By (\ref{eqn-annihilator}) the ordinary homology is
the submodule $H_*(X_{\gamma})= H^{T(\mathbb C)}_*(X_{\gamma})\{ \mathcal I\}$
which is killed by the augmentation ideal $\mathcal I,$ so {\it we obtain an
action of $(\widetilde W \rtimes \mathbf{Aut})_{\gamma}$ on the ordinary
homology of the affine Springer fiber.}

We remark that in many cases (including $G = SL(n)$, $Sp(n)$, or $O(n)$) it is
possible to construct a (continuous) action of the group
(\ref{eqn-automorphismgroup}) on $X$ which is $T(\mathbb C)$ covariant,
meaning that $\tau(t\cdot x) =\tau(t)\cdot \tau(x)$ for $t\in T(\mathbb C)$,
$x \in X$, and $\tau \in \widetilde{W} \rtimes\mathbf{Aut}.$  The induced
action on $H_*^{T(\mathbb C)}(X)$ then agrees with the left action.  However,
even in these cases, the topological action of the subgroup (\ref{eqn-Wgamma})
does not necessarily preserve the Springer fiber $X_{\gamma}.$

\section{Endoscopic groups}\label{sec-endoscopic}
\subsection{}\label{subsec-endoscopic1}
In this section we assume, for simplicity only, that $G$ is adjoint.  Fix
endoscopic data $(H,s)$ for $G$.  This means that $H$ is a connected reductive
complex algebraic group, that $s \in \widehat{T} \subset \widehat{G}$, and
$\widehat{H} = \widehat{G}_s$ is the centralizer of $s$ in $\widehat{G}.$
Then $\widehat{H}$ is connected since $\widehat{G}$ is simply connected
(\cite{Humphreys} \S 2.11).  We assume moreover that Borel pairs $T \subset B
\subset G$ and $T_H \subset B_H \subset H$ (defined over $\mathbb C$) have
been chosen, giving rise to a canonical isomorphism
\begin{equation}\label{eqn-torusiso}
  T_H \cong T.\end{equation}
Then $\widehat{T} = T_{\widehat{H}}$ is also a maximal torus in $\widehat{H}.$
The set of co-roots of $T$ in $H$ is
\begin{equation*}
\Phi^{\vee}(H,T_H) = \left\{ \alpha^{\vee} \in \Phi^{\vee}(G,T) :\
s(\alpha^{\vee}) = 1 \right\}.
\end{equation*}
Let $\mathcal D = \mathcal D(\mathfrak t)\cong \bS(\mathfrak t^*)$ be the
graded algebra of differential operators on $\mathfrak t^*.$  Each $\alpha \in
\Phi(G,T)$ determines a differential operator $\partial_{\alpha}\in \mathcal
D$ of degree 1.  Define the {\it homological transfer factor}
\begin{equation}\label{eqn-transfer}
  \Delta = \prod_{\alpha \in \Phi^+(G) - \Phi^+(H)} \partial_{\alpha}^{\val
\alpha'(\gamma)} \in \mathcal D = H^*_{T(\mathbb C)}(\text{pt}).\end{equation}
It is a homogeneous element of some degree, $r=\sum_{\alpha \in \Phi^+(G) -
\Phi^+(H)}\val \alpha'(\gamma).$  The action of $\Delta$ is a surjection of
graded $\mathcal D$ modules
\begin{equation*}\begin{CD}
\bS(\mathfrak t) @>{\Delta \cdot}>> \bS(\mathfrak t)[-r] @>>> 0
\end{CD}\end{equation*}
(where $A[-r]$ denotes the shift in grading, $A[-r]_k = A_{k-r}$ for any
graded vector space $A = \bigoplus_k A_k$).

Fix $\gamma \in \mathfrak t(\mathfrak o)$ regular and semisimple and let
$\gamma_H \in \mathfrak t_H = \text{Lie}(T_H)$ be the element which
corresponds to $\gamma$ under the isomorphism (\ref{eqn-torusiso}).  Then
$\gamma_H$ acts as a ``vector field'' on the affine Grassmannian $X^H = H(F)
/H(\mathfrak o)$ and we denote its fixed point set by $X^H_{\gamma_H}.$

Let $\mathbf{Aut}^{G,H}$ be the subgroup of the automorphism group of
the based root
datum for $G$ which preserves the roots of $H.$  Let $W^H$ be the Weyl
group for $H,$ set $\widetilde{W}^H = \Lambda \rtimes W^H$ and
$\widetilde{\mathcal W}^{G,H}=
\widetilde{W}^H \rtimes \mathbf{Aut}^{G,H}$, and define a subgroup of
$\widetilde{\mathcal
W}^{G,H}$ by
\begin{equation} \label{eqn-WGHgamma}
  \widetilde{\mathcal W}^{G,H}_{\gamma} =
\left\{ \tau \in  \widetilde{\mathcal W}^{G,H}: \
\val (\tau \alpha(\gamma)) = \val(\alpha(\gamma)) \text{ for all } \alpha \in
\Phi^+(G,T)  \right\}.  \end{equation}
Note that the square of $\Delta$ is invariant under
$\widetilde{\mathcal W}^{G,H}_{\gamma}$ and hence that there exists a
sign character
\[
\eta:\widetilde{\mathcal W}^{G,H}_{\gamma} \to \{\pm 1\}
\]
such that $\tau(\Delta)=\eta(\tau) \cdot \Delta$ for all $\tau \in
\widetilde{\mathcal
W}^{G,H}_{\gamma}$.

Let $J \subset \mathbb C[\Lambda]$ be the multiplicative subset which is
generated by the elements $1 - \alpha^{\vee}$ for $\alpha^{\vee} \in
\Phi^{\vee}(G,T) - \Phi^{\vee}(H,T_H).$ For any $\mathbb C[\Lambda]$ module
$M$, let $ M_J = J^{-1}M $ denote the $\mathbb C[\Lambda]_J$ $=$ $J^{-1}
\mathbb C[\Lambda]$ module which is obtained by inverting the elements of $J.$
The localization of the module $\mathbb C[\Lambda] \otimes \bS(\mathfrak t)$
is $ \mathbb C[\Lambda]_J \otimes \bS(\mathfrak t).$

\begin{thm}\label{thm-mainA} Suppose the homology of $X_{\gamma}$ and of
$X^H_{\gamma_H}$ is pure. Then the surjection
\begin{equation}\label{eqn-otimesdelta} 1 \otimes \Delta: \mathbb C[\Lambda]
\otimes_{\mathbb C} \bS(\mathfrak t) \to \mathbb C[\Lambda]\otimes_{\mathbb C}
\bS(\mathfrak t)[-r] \end{equation}
induces a
homomorphism of $\mathcal D$ modules \begin{equation}\label{eqn-Psi}  \Psi:
H_*^{T(\mathbb C)}(X_{\gamma}) \to H_*^{T(\mathbb C)}(X^H_{\gamma_H})[-2r]
\end{equation}
which becomes an isomorphism
\[ \Psi_J: H_*^{T(\mathbb C)}(X_{\gamma})_J \to H_*^{T(\mathbb C)}
(X^H_{\gamma_H})_J[-2r]  \]
after localizing with respect to $J.$
The homomorphisms $\Psi$ and $\Psi_J$ are not quite
$\widetilde{\mathcal W}^{G,H}_{\gamma}$
equivariant; rather they transform under $\widetilde{\mathcal
W}^{G,H}_{\gamma}$ by the
sign character $\eta$ defined above.

  \end{thm}

Assuming the homology of $X_{\gamma}$ is pure, by Lemma \ref{lem-pureind} and
equation (\ref{eqn-annihilator}) the ordinary homology of $X_{\gamma}$ is
given by the $\mathcal D$-submodule of the equivariant homology
\[ H_*(X_{\gamma}; \mathbb C) \cong H_*^{T(\mathbb C)}(X_{\gamma}; \mathbb C)
\left\{ \mathcal I \right\}\]
which is annihilated by the augmentation ideal $\mathcal I \subset \mathcal
D.$  We conclude:

\begin{cor}\label{cor-ordinaryH}
The action of $1\otimes\Delta$ on the equivariant homology induces a
homomorphism $H_*(X_{\gamma}) \to H_*(X^H_{\gamma_H})[-2r]$ and an isomorphism
\begin{equation}\label{eqn-ordinaryH}
  H_*(X_{\gamma}; \mathbb C)_J \cong H_*(X^H_{\gamma_H}; \mathbb C)_J[-2r].
\end{equation}
These maps again transform under $\widetilde{\mathcal W}^{G,H}_{\gamma}$ by the
sign character $\eta$.
\end{cor}
Although the homomorphism (\ref{eqn-Psi}) may be described relatively easily
in equivariant homology, the resulting isomorphism (\ref{eqn-ordinaryH}) in
ordinary homology is much more complicated.

\subsection{Proof of Theorem \ref{thm-mainA}}
The mapping $\Psi$ is well defined because the mapping
(\ref{eqn-otimesdelta}) kills the submodule
\[ L_{\alpha,\gamma} = \sum_{d=1}^{\val(\alpha'(\gamma))} (1-\alpha^{\vee})^d
\mathbb C[\Lambda] \otimes \bS(\mathfrak t)\left\{\partial_{\alpha}^d\right\}
\]
whenever $\alpha \in \Phi^+(G) - \Phi^+(H).$  If $w \in \widetilde{W} \rtimes
\mathbf{Aut}$ and if $f \in \mathbf{S}(\mathfrak t)$ then $w(\Delta f) =
w(\Delta) w(f).$  If $w$ lies in the subgroup
$\widetilde{\mathcal W}^{G,H}_{\gamma}$  then $w(\Delta) = \eta(w)
\cdot \Delta.$  Therefore
the maps $1 \otimes \Delta$ and $\Psi$
transform under $\widetilde{\mathcal W}^{G,H}_{\gamma}$ by the
sign character $\eta$.

Now let us check that $\Psi_J$ is an isomorphism.  The surjection
\[ (1\otimes \Delta)_J: \mathbb C[\Lambda]_J \otimes \bS(\mathfrak t)
\to  \mathbb C[\Lambda]_J \otimes \bS(\mathfrak t)[-r] \]
has kernel
\[ \ker(1\otimes \Delta)_J = \mathbb C[\Lambda]_J \otimes \bS(\mathfrak t)
\left\{ \Delta
\right\} \]
and it restricts to a mapping $h_J$ in the following diagram.
\begin{equation*}\begin{CD}
0 &\longrightarrow \ & \sum_{\alpha\in\Phi^+(G)}J^{-1} L_{\alpha,\gamma}
&\longrightarrow \ &  \mathbb C[\Lambda]_J \otimes \bS(\mathfrak t)
&\longrightarrow \ & H^{T(\mathbb C)}_*(X_{\gamma})_J
& \longrightarrow \ & 0 \\
&& @V{h_J}VV  @VV{(1\otimes\Delta)_J}V &\\
0 &\longrightarrow \ & \sum_{\alpha \in \Phi^+(H)} J^{-1} L_{\alpha,\gamma}
[-r] &\longrightarrow \  & \mathbb C[\Lambda]_J \otimes \bS(\mathfrak t)[-r]
&\longrightarrow \ & H^{T(\mathbb C)}_*(X^H_{\gamma_H})_J[-2r]
&\longrightarrow \ & 0
\end{CD}\end{equation*}
It suffices (by the snake lemma) to show that $h_J$ is surjective, and that
\begin{equation}\label{eqn-kerDelta} \ker(1\otimes\Delta)_J \subset
\sum_{\alpha \in \Phi^+(G)} J^{-1} L_{\alpha,\gamma}.
\end{equation}
If $\alpha \in \Phi^+(G) - \Phi^+(H)$ (that is, if $s(\alpha^{\vee}) \ne 1$)
then $(1-\alpha^{\vee})$ becomes invertible in $\mathbb C[\Lambda]_J.$  If $a
\le b$  then $\bS (\mathfrak t) \left\{\partial_{\alpha}^a\right\} \subset
\bS(\mathfrak t) \left\{\partial_{\alpha}^{b} \right\}.$  Therefore the
localization of the submodule $ L_{\alpha,\gamma}$ is
\[
J^{-1}L_{\alpha,\gamma} = \mathbb C[\Lambda]_J \otimes \bS(\mathfrak t)
\left\{\partial_{\alpha}^{\val(\alpha'(\gamma))}\right\}.
\]
By Lemma \ref{lem-relativelyprime}, the sum
\[ \sum_{\substack {{\alpha}\in \Phi^+(G)\\ s(\alpha^{\vee}) \ne 1}}
\bS(\mathfrak t) \left\{\partial_{\alpha} ^{\val(\alpha'(\gamma))}\right\} =
\bS(\mathfrak t) \left\{\Delta\right\}
\]
is killed by $\Delta.$  Therefore
\[ \ker(1\otimes\Delta)_J = \sum_{\substack {{\alpha}\in \Phi^+(G)\\
s(\alpha^{\vee}) \ne 1}} J^{-1} L_{\alpha,\gamma}\]
which proves (\ref{eqn-kerDelta}).  On the other hand, if $\alpha \in
\Phi^+(G)$ and $s(\alpha^{\vee})=1$ then
\[ J^{-1}L_{\alpha,\gamma} = \sum_{d=1}^{\val(\alpha'(\gamma))}
(1-\alpha^{\vee})^d \mathbb C[\Lambda]_J \otimes \bS(\mathfrak t)
\left\{\partial_{\alpha}^d\right\}.\]
In this case (by Lemma \ref{lem-relativelyprime}) the action of $\Delta$ is a
surjection
\[\bS(\mathfrak t)\left\{\partial_{\alpha}^d\right\} \to
\bS(\mathfrak t)\left\{\partial_{\alpha}^d\right\}[-r]\]
which implies that $h_J$ is surjective.  \qed

\subsection{Remark} We did not use the full strength of the endoscopic
hypothesis on~$H$, and in fact Theorem \ref{thm-mainA} remains valid whenever
$H$ is a connected reductive group for which there exists a Borel pair $T_H
\subset B_H$ and an isomorphism $T \cong T_H$ such that the set of
positive roots $\Phi^+(H,T_H)$ is a subset of the positive roots $\Phi^+(G,T)$.

\section{The quotient under $\Lambda$}\label{sec-quotient}
\subsection{}
As in \S \ref{sec-endoscopic}, fix a Borel pair $T \subset B \subset G,$ with
resulting positive roots $\Phi^+.$ Fix regular
$\gamma \in \mathfrak t(\mathfrak o)$
with affine Springer fiber $X_{\gamma}$ and its lattice of translations
$\Lambda.$  The quotient $\Lambda \backslash X_{\gamma}$ admits the structure
of a complex projective algebraic variety \cite{KL}.  Let $(H,s)$ be
endoscopic data, with $s \in \widehat{T}$, and assume that $s$ has finite order.
  View $s$ as a character $s:\Lambda
\to \text{GL}_1(\mathbb C)$ and let $\mathbb C_s$ be the 1-dimensional
$\mathbb C[\Lambda]$ module defined by $s.$  Then $s$ determines 1-dimensional
local systems
\[ \mathcal L_s = X_{\gamma} \times_{\Lambda} \mathbb C_s\ \text{ and }
\mathcal M_s = X^H_{\gamma_H} \times_{\Lambda} \mathbb C_s\]
on $\Lambda \backslash X_{\gamma}$ and $\Lambda \backslash X^H_{\gamma_H}$
respectively.  Let $r=\deg(\Delta)$ as in (\ref{eqn-transfer}).

\begin{prop}\label{prop-E2} Suppose that $H_*(X_{\gamma}; \mathbb C)$ and
$H_*(X^H_{\gamma_H}; \mathbb C)$ are pure. Then there is an isomorphism
\begin{equation}\label{eqn-equalquotients}
H_m(\Lambda \backslash X_{\gamma}; \mathcal L_s) \cong H_{m-2r}
(\Lambda \backslash X^H_{\gamma_H}; \mathcal M_s).  \end{equation}
If the homology of $X_{\gamma}$ is both pure and Tate, and if the same is true
of $X^H_{\gamma_H}$, then this isomorphism may be chosen in a canonical manner.
\end{prop}

\subsection{Proof}
   The $E^2$ term of the Cartan-Leray spectral sequence for the above
$\Lambda$-equivariant homology group is
\[ E^2_{pq}(X_{\gamma};\mathbb C_s) =
H_p(\Lambda; H_q(X_{\gamma}; \mathbb C_s)) = \text{Tor}_p^{\mathbb C[\Lambda]}
(H_q(X_{\gamma}), \mathbb C_s) \implies H_{p+q}(\Lambda \backslash X_{\gamma};
\mathcal L_s).\]
By Lemma \ref{lem-collapse} the spectral sequence collapses and gives rise to
an isomorphism
\begin{equation}\label{eqn-Tor}
H_m(\Lambda \backslash X_{\gamma}; \mathcal L_s) \cong
\bigoplus_{p+q=m} \text{Tor}_p^{\mathbb C[\Lambda]}(H_q(X_{\gamma};\mathbb C),
\mathbb C_s)
\end{equation}
which is canonical if $H_*(X_{\gamma};\mathbb C)$ is also Tate.

Let $J \subset \mathbb C[\Lambda]$ be the multiplicative subset generated by
the collection of elements $(1 - \alpha^{\vee})$ with $\alpha \in \Phi^+(G)$
and $s(\alpha^{\vee}) \ne 1.$  The homomorphism $\mathbb C[\Lambda] \to
\mathbb C_s$ ($\ell \mapsto s(\ell)$) which gives rise to the local system
$\mathcal L_s$ may be factored as the composition of ring homomorphisms
\begin{equation}\label{eqn-factoring}
\mathbb C[\Lambda] \to  C[\Lambda]_J \to \mathbb C_s \end{equation}
in which the second map is defined by
\[ (1-\alpha^{\vee})^{-1} \ell \mapsto (1-s(\alpha^{\vee}))^{-1} s(\ell) \]
for any $\alpha^{\vee} \in \Phi^{\vee}(G,T) - \Phi^{\vee}(H, T_H)$ and any
$\ell \in \Lambda.$ It is well defined because $s(\alpha^{\vee}) \ne 1.$  For
any $\mathbb C[\Lambda]$ module $V$, the composition (\ref{eqn-factoring})
induces an isomorphism
\[
V \otimes_{\mathbb C[\Lambda]} \mathbb C_s \cong
V \otimes_{\mathbb C[\Lambda]} \mathbb C[\Lambda]_J \otimes_{\mathbb
C[\Lambda]_J} \mathbb C_s \cong
V_J \otimes_{\mathbb C[\Lambda]_J} \mathbb C_s.\]
Taking derived functors gives an isomorphism
$\text{Tor}^{\mathbb C[\Lambda]}_p(V,\mathbb C_s) \cong \text{Tor}^{\mathbb
C[\Lambda]_J}_p(V_J, \mathbb C_s).$  Taking $V = H_q(X_{\gamma}; \mathbb C)$
gives an
isomorphism
\begin{equation*}
H_m(\Lambda \backslash X_{\gamma}; \mathcal L_s) \cong \bigoplus_{p+q=m}
\text{Tor}_p^{\mathbb C[\Lambda]_J}
(H_q(X_{\gamma}; \mathbb C)_J, \mathbb C_s).  \end{equation*}
So  Proposition \ref{prop-E2} follows from Corollary \ref{cor-ordinaryH}.
\qed

\subsection{}
Although we have not constructed the isomorphism (\ref{eqn-equalquotients}) in
a canonical way,
the isomorphism induced by $\Delta,$
\begin{equation*}
  E^2_{pq}(X_{\gamma}; \mathbb C_s) \cong
E^2_{p,q-2r}(X^H_{\gamma_H}; \mathbb C_s) \end{equation*}
is canonical, and it transforms under $\widetilde{\mathcal W}^{G,H}_{\gamma}$
by the
sign character $\eta$.  It follows that the Lefschetz number of any $w\in
\widetilde{\mathcal
W}^{G,H}_{\gamma}$ on $E^2_*(X_{\gamma};\mathbb C_s)$ is equal to that on
$E^2_*(X^H_{\gamma_H}; \mathbb C_s)$ times the sign $\eta(w).$

\section{Combinatorial Lemmas}\label{sec-combinatorial}
\subsection{}
Let $\Lambda$ be a free abelian group of rank one, and let $\mathbb Q[\Lambda]
$ be its group ring.  (The group law in $\Lambda$ will be written
multiplicatively.)  A choice of generator $\alpha^{\vee} \in \Lambda$
determines an isomorphism $\Lambda \cong \mathbb Z.$  Write $\ell_a \in
\Lambda$ for the element corresponding to $a \in \mathbb Z.$  Then
multiplication by $\ell_1 = \alpha^{\vee}$  acts as a ``shift operator'':
$\alpha^{\vee} \ell_a = \ell_{a+1}.$

The algebra $\mathbb Q[\partial_x,\partial_t]$ of differential operators in
two variables acts on the algebra $\mathbb Q[x,t]$ of polynomials in two
variables, and hence also on the algebra $\mathbb Q[\Lambda] \otimes_{\mathbb
Q}\mathbb Q[x,t]$ with
\[ \ker(\partial_t) = \mathbb Q[\Lambda] \otimes_{\mathbb Q} \mathbb Q[x].\]
Fix $d,m \in \mathbb Z$ with $d \ge 1.$  Define
\begin{equation}\label{eqn-Fsum} f_{m,d} =
\sum_{m \leq a < b \leq m+d} C_{ab}(\ell_b - \ell_a) \otimes \left( (a+b)x-t)
\right)^{d-1}
\in \mathbb Q[\Lambda] \otimes_{\mathbb Q} \mathbb Q[x,t]\end{equation}
where $C_{ab} \in \mathbb Z$ is the integer
\[
C_{ab} = (-1)^{a-b}(a-b) {\binom{d}{a-m} \binom{d}{b-m}}.\]
Define $J_v\subset \mathbb Q[\Lambda] \otimes_{\mathbb Q} \mathbb Q[x,t]$ to
be the span
\[ J_v = \sum_{d=1}^v \sum_{m\in \mathbb Z} \mathbb Q f_{m,d}. \]

\begin{lem}\label{lem-Fmd}  Fix $d,m \in \mathbb Z$ with $d \ge 1.$
Then
\begin{equation}\label{eqn-Fmd}
f_{m,d} = (-1)^dd! (1-\alpha^{\vee})^d \ell_{m} \otimes x^{d-1}\in
\ker(\partial_t).
\end{equation}
If $v \ge 1$ then
\begin{equation}\label{eqn-module}
J_v = \sum_{d=1}^v
(1-\alpha^{\vee})^d\mathbb Q[\Lambda] \otimes_{\mathbb Q} \mathbb Q[x]\left\{
\partial_x^d \right\} \end{equation}
where $\mathbb Q[x]\left\{ \partial_x^d \right\}$ denotes the polynomials
which are annihilated by $\partial_x^d$, that is, the polynomials of degree
$\le d-1.$
\end{lem}

\subsection{Proof}
We use the fact \cite{Hall} \S 1.1 that for any polynomial $p$,
\begin{equation}\label{eqn-binom}
\sum_{k=0}^n (-1)^k \binom{n}{k} p(k) = \begin{cases}
0 & \text{if }\deg{p} \le n-1 \\
(-1)^nn! & \text{if } p(k) = k^n. \end{cases} \end{equation}
Since $C_{ab} = - C_{ba}$ we may express $f_{m,d}$ as:
\begin{align*}
&\sum_{a=m}^{m+d} \sum_{b=m}^{m+d} C_{ab} \ell_b \otimes \bigl( (a+b)x -
t\bigr)^{d-1}\\
= &\sum_{a=m}^{m+d} \sum_{b=m}^{m+d} \left( \ell_b \otimes C_{ab} \sum_{j=0}
^{d-1} \binom{d-1}{j} (a+b)^j x^j (-1)^{d-1-j} t^{d-1-j} \right)\\
=&\sum_{b=m}^{m+d} \sum_{j=0}^{d-1} \left( \ell_b \otimes (-1)^{b'} \binom{d}
{b-m}\binom{d-1}{j} t^{d'} \sum_{a=m}^{m+d} (-1)^{a-m} \binom{d}{a-m} (a-b)
(a+b)^j x^j \right)
\end{align*}
where $b' = b+d-1-j-m$ and $d'=d-1-j.$  By (\ref{eqn-binom}), the innermost
sum vanishes unless $j=d-1$, leaving
\[ f_{m,d} = (-1)^d d! \sum_{b=m}^{m+d}\left((-1)^{b-m} \binom{d}{b-m} \ell_b
\right) \otimes x^{d-1}\]
which is (\ref{eqn-Fmd}).   To verify the reverse inclusion in equation
(\ref{eqn-module}) it suffices to show that $(1-\alpha^{\vee})^d \ell_{m}
\otimes x^{e-1}$ is in $J_v$ for every $d,e$ with $1 \le e \le d \le v.$ But
  \begin{align*}
(1-\alpha^{\vee})^d \ell_{m} \otimes x^{e-1} &= (1-\alpha^{\vee})
^e(1-\alpha^{\vee})^{d-e}
\ell_{m} \otimes x^{e-1} \\
&= (1-\alpha^{\vee})^e \sum_{j=0}^{d-e}(-1)^{d-e-j}\binom{d-e}{j}\ell_{m+j}
\otimes x^{e-1}\\
&= \frac{(-1)^{d-e}}{e!} \sum_{j=0}^{d-e} (-1)^j \binom{d-e}{j} f_{m+j,e}
\end{align*}
by (\ref{eqn-Fmd}).  Since $e \le v$, this quantity lies in $J_v.$  This
completes the proof. \qed

\begin{lem}\label{lem-degree}
Fix $d,h,v \ge 1$ with $v \le h.$  Fix $m \in \mathbb Z.$  Let
\begin{equation}\label{eqn-g} g = \sum_{\substack{m \le a < b \le m+h\\ b-a
\le v}}(\ell_b - \ell_a) \otimes G_{ab}\left( (a+b)x-t\right)^{d-1}
\end{equation}
where $G_{ab} \in \mathbb Q$ \parens{with $b > a$}.  If $\partial_tg =0$ and
$d > v$  then $g= 0.$
\end{lem}

\subsection{Proof}  Write $g_{ab} = G_{ab}\left( (a+b)x-t\right)^{d-1}.$
The sum in (\ref{eqn-g}) may be written either as
\[ \sum_{a=m}^{m+h-1}\ \sum_{b=a+1}^{\min(m+h,a+v)} \text{ or }
\sum_{b=m+1}^{m+h}\ \sum_{a=\max(m,b-v)}^{b-1}.\]
  Therefore
\begin{align*}
g &= \sum_{b=m+1}^{m+h}\  \sum_{a=\max(m,b-v)}^{b-1} \ell_b \otimes g_{ab}
  -\sum_{a=m}^{m+h-1}\ \sum_{b=a+1}^{\min(m+h,a+v)}\ell_a \otimes g_{ab} \\
&= \sum_{a=m}^{m+h} \ell_a \otimes \left(
\sum_{b=\max(m,a-v)}^{a-1} g_{ba} - \sum_{b=a+1}^{\min(m+h,a+v)} g_{ab}\right)
\end{align*}
(but if $a=m$ then the first sum is empty, while if $a = m+h$ the second sum
is empty). Fix $a$ with $m \le a \le m+h.$  Expanding $((a+b)x-t)^{d-1},$ the
equation $\partial_tg=0$ gives the following system of linear equations:
\begin{equation}\label{eqn-system}
\sum_{b=\max(m,a-v)}^{a-1} G_{ba}(a+b)^j - \sum_{b=a+1}^{\min(m+h,a+v)}G_{ab}
(a+b)^j =
0\end{equation}
for $ 0 \le j \le d-2.$  This is a system of $d-1$ homogeneous equations of
van der Monde type, in $\min(m+h,a+v)- \max(m,a-v)$ variables. If the number
of equations equals or exceeds the number of variables, then only the zero
solution exists, since the corresponding van der Monde determinant is nonzero.
Suppose $d > v.$ First consider the system (\ref{eqn-system}) corresponding to
$a = m.$ Then the first sum in (\ref{eqn-system}) is empty, leaving the second
sum, which consists of $d-1$ (homogeneous) equations (of van der Monde type)
in $\le v$ unknowns $G_{mb}$ (with $a+b \le b \le a+v$).  So these
coefficients vanish.  Suppose by induction that the coefficients $G_{ab}$ have
been shown to vanish for all pairs $(a,b)$ with $a <t$ and $a<b\le
\min(m+h,a+v).$  Consider the system (\ref{eqn-system}) in the case $a = t.$
The terms in the first sum, $\sum_{x=\max(m,t-v)}^{t-1}G_{xt}$ vanish by the
induction hypothesis, leaving only the second sum.  This is a system of $d-1$
homogeneous equations in $\le v$ variables so only the zero solution $G_{tb}=
0$ exists.  This completes the induction, so we conclude that $g=0.$  \qed

We remark that the same (inductive) argument shows that if $d=v=h$ then there
is at most a one dimensional space of solutions, so $g$ is a multiple of
$f_{m,d}.$

\subsection{}\label{subsec-Pv}
Define $P_v \subset \mathbb Q[\Lambda] \otimes_{\mathbb Q} \mathbb Q[x,t]$ to
be the vector space spanned by elements
\begin{equation}\label{eqn-gab}
(\ell_b - \ell_a) \otimes g_{ab}((a+b)x -t)\end{equation}
where $g_{ab}$ are polynomials, and where $|b-a| \le v.$  Let $P_v \left\{
\partial_t \right\} = P_v \cap \ker(\partial_t).$

\begin{prop} \label{prop-Fmdspan} Fix $v \ge 1.$  Then
$P_v\left\{\partial_t\right\} = J_v,$
that is,
\[\ker(\partial_t)\cap \sum_{|b-a| \le v} \mathbb Q(\ell_b - \ell_a) \otimes
\mathbb Q[(a+b)x- t]  = \sum_{d =1}^v (1-\alpha^{\vee})^d \mathbb Q[\Lambda]
\otimes \mathbb Q[x] \left\{ \partial_x^d \right\}. \]
\end{prop}

\subsection{Proof}
The inclusion $J_v \subset P_v\left\{ \partial_t \right\}$ follows from Lemma
\ref{lem-Fmd}, so we only need to verify the reverse inclusion.  Let $P_{v,h}$
be the vector subspace of $P_v$ spanned by elements (\ref{eqn-gab}) such that
each $g_{ab}(z)= G_{ab}z^h$ is homogeneous of degree $h$ (with $G_{ab} \in
\mathbb Q$).  Set $P_{v,h}\left\{ \partial_t \right\} = P_{v,h} \cap \ker
\partial_t.$  Then
\[P_v\left\{ \partial_t \right\} = \sum_{h \ge 0} P_{v,h} \left\{ \partial_t
\right\};\]
and Lemma \ref{lem-degree} says that $P_{v,h}\left\{ \partial_t \right\} = 0$
whenever $h \ge v.$  So we need to show that $P_{v,d-1}\left\{ \partial_t
\right\} \subset J_v$ whenever $d \le v.$  Since both $P_v\left\{ \partial_t
\right\}$ and $J_v$ are modules over $\mathbb Q[\partial_x]$ it suffices to
consider the case $d=v,$ that is, we must show that $P_{v,v-1} \left\{
\partial_t \right\} \subset J_v.$  Let
\begin{equation}\label{eqn-gd} g = \sum_{\substack{m \le a < b \le m+N \\ b -
a \le v}} (\ell_b - \ell_a) \otimes G_{ab} ((a+b)x-t)^{v-1}\in P_{v,v-1}
\left\{ \partial_t \right\}
\end{equation}
(for some $m,N\in\mathbb Z$) where $G_{ab} \in \mathbb Q.$  The expression
(\ref{eqn-Fsum}) for $f_{m,v}$ contains a unique term with $b-a=v.$   However
the sum (\ref{eqn-gd}) for $g$ contains only terms with $b-a \le v$, so by
subtracting appropriate multiples of $f_{*,v}\in J_v \cap P_{v,v-1}
\{\partial_t\}$ we can eliminate all terms with $b-a =v.$  In other words,
there exists $h \in J_v \cap P_{v,v-1} \left\{ \partial_t \right\}$ such that
$g - h \in P_{v-1,v-1}\left\{ \partial_t \right\} = 0.$ Therefore $g \in J_v.$
\qed

\subsection{Remarks}\label{subsec-constellations}
The lemmas in this section refer to the equivariant homology $H_*^{T(\mathbb
C)}(X^{\SL}_{\le v})$ of the affine Springer fiber.  Each 1-dimensional orbit
$E_{ab} \subset X^{\SL}_{\le v}$ of the extended torus $\widetilde{T}(\mathbb
C)$ determines a relation in $\mathbb C[\Lambda] \otimes \mathbb C[x,t] =
H_*^{\widetilde{T}(\mathbb C)}(\Lambda).$  For $1 \le d \le v$ the sum
(\ref{eqn-Fsum}) corresponds to a ``constellation'' of (finitely many)
1-dimensional orbits, such that the resulting relation $f_{m,d}$ lies in the
subalgebra $\mathbb C[\Lambda] \otimes  \mathbb C[x] = H_*^{T(\mathbb C)}
(\Lambda).$  This constellation consists of all the 1-dimensional orbits which
are contained in a single irreducible component of $X^{\SL}_{\le d}.$  One may
think of this constellation of 1-dimensional orbits as being attached to the
single orbit (the ``longest one'' in the constellation) which joins the
lattice points $\ell_m$ and $\ell_{m+d}.$  So $H_*^{T(\mathbb C)}(X_{\le v})$
is the quotient of $\mathbb C[\Lambda] \otimes \mathbb C[x]$ by the relations
(\ref{eqn-Fmd}), one for each 1-dimensional $\widetilde{T}(\mathbb C)$-orbit
in $X_{\le v}.$

\section{Affine flag manifold for $\SL$}\label{sec-flagSL2}
   The results in this section are parallel to those of \S \ref{sec-SL2},
\S \ref{sec-EqHoSL2}, \S \ref{sec-ssrankone}, and \S \ref{sec-combinatorial}
so they will be presented without detailed proofs.

\subsection{} Throughout this section we take $G = \SL.$  Let $B =
\left(\begin{smallmatrix} * & * \\ 0 & *  \end{smallmatrix}\right)$ be the
standard Borel subgroup of $\SL$ with corresponding Iwahori
subgroup $ I \subset G(\mathfrak o)$ and let $Y=Y^{\SL} = G(F)/  I$ be the
affine flag manifold with basepoint $x_0 = I.$ Let $T \subset G$ be the torus
of diagonal matrices, with its lattice of translations $\Lambda =\left\{
\alpha^{\vee}(\epsilon^m):\ m \in \mathbb Z \right\}$ where $\alpha^{\vee} \in
\Phi^{\vee}$ is the simple co-root determined by $T \subset B$ (cf. \S
\ref{sec-affinespringer}).  The torus $T(\mathbb C)$ acts on $Y$ with fixed
points
\begin{equation}\label{eqn-morefixedpoints}
\ell_m = \left( \begin{matrix} \epsilon^m & 0 \\ 0 & \epsilon^{-m}
\end{matrix}
\right) I \ \text{ and }\ r_m = \left( \begin{matrix} 0 & \epsilon^m \\
-\epsilon^{-m} & 0 \end{matrix} \right) I
\end{equation}
for $m \in \mathbb Z.$  Let $W=W(G,T) = \left\{ 1,w_{\alpha}\right\}$ be the
Weyl group of $G.$ Since the extended affine Weyl group $\widetilde{W}=
\Lambda \rtimes W$ acts simply transitively on these fixed points, the choice
of basepoint $x_0 = \ell_0 \in Y$ determines an identification of
$\widetilde{W}$ with this set of fixed points, which identifies $\alpha^{\vee}
(\epsilon^m) \in \Lambda$ with $\ell_m \in Y,$ and identifies the reflection
$w_{\alpha}$ with $r_0 \in Y.$ Then  $\ell_a \ell_b = \ell_{a+b}$, $\ell_a r_b
= r_{a+b}$,  $r_a \ell_b = r_{a-b},$ and $r_ar_b = \ell_{a-b}.$

For $m \le 0$ define
\begin{equation}\label{eqn-xmym}
x_m = \left( \begin{matrix} 1 & \epsilon^m \\ 0 & 1 \end{matrix} \right)
  I \  \text{ and }\
y_m = \left( \begin{matrix} \epsilon^m & -1 \\ 1 & 0 \end{matrix} \right)  I.
\end{equation}
Then $x_0 = \ell_0$ is still the basepoint.  For notational brevity, put $y_1
= r_0.$ Let $\widetilde{T}(\mathbb C)$ denote the extended torus
(\ref{eqn-extendedtorus}). If $m \le -1$ and $s \in \mathbb Z,$ the
$\widetilde{T} (\mathbb C)$ orbit of $\ell_sx_m$ is 1-dimensional, and it
connects the fixed points $\ell_s$ and $r_{m+s}.$  If $m \le 0,$ the
$\widetilde{T}(\mathbb C)$ orbit of $\ell_s y_m$ is 1-dimensional and it
connects the fixed points $r_s$ and $\ell_{m+s}.$  Let $\pi:Y^{\SL} \to
X^{\SL}$ be the natural projection from the affine flag manifold to the affine
Grassmannian. It is a $G(F)$-equivariant fibration with fiber isomorphic to
$G(\mathbb C) / B(\mathbb C) \cong \mathbb C \mathbb P^1.$  It satisfies
$\pi(\ell_m) = \pi(r_m)$ for all $m,$ and $\pi(x_m) = \pi(y_m)$ for all $m\le
0.$

For $m \le 0$ the $T(F)$-orbit of $x_m\in Y$ is $(-m)$ dimensional and it
projects isomorphically under $\pi$ to the $T(F)$-orbit of the point $\pi(x_m)
$ (which was denoted $x_m$ in \S \ref{sec-SL2}).  For $m \le 1$ the $T(F)$
orbit of $y_m\in Y$ is $1-m$ dimensional.  If $m \le 0,$ it fibers over the
$T(F)$ orbit of $\pi(y_m) = \pi(x_m)$ with fiber a 1-dimensional affine space.

Fix a regular element $\gamma \in \mathfrak t(\mathfrak o)$ and let
\[ Y_{\gamma} = \left\{ x  I \in Y:\ \text{Ad}(x^{-1})(\gamma) \in \text{Lie}
(I) \right\}
\]
be the affine Springer fiber corresponding to $\gamma$ in the affine flag
manifold $Y.$  It is preserved by the extended torus $\widetilde T(\mathbb
C),$ and the mapping $\pi:Y_{\gamma} \to X_{\gamma}$ is surjective (however it
may fail to be a fibration).  The homology of $Y_{\gamma}$ is pure
(\cite{GKM2}). The paving of $Y$ by affine cells gives rise to a
paving of $Y_{\gamma}$
by affine cells. In fact, the intersection of $Y_{\gamma}$ with a
single Bruhat cell in $G(F)/I$
is either empty, or else it is a single affine cell.

\begin{prop}
The affine flag manifold $Y$ is the disjoint union of countably many $T(F)$
orbits,
\begin{equation}\label{eqn-moretorusorbits}
Y = \coprod_{n \le 0} T(F)\cdot \left\{x_n \cup y_{n+1} \right\}.
\end{equation}
The affine Springer fiber $Y_{\gamma}$ is the union
\begin{equation*}
Y_{\gamma} = Y_{\le v}= \bigcup_{n = -v}^{0} T(F) \cdot \left\{ x_n \cup
y_{n+1} \right\} \end{equation*}
where $v = \val(\alpha'(\gamma)).$  For each ordered pair of integers $s,t \in
\mathbb Z$ there is a unique one dimensional orbit $E_{st}\subset Y$ of the
extended torus  $\widetilde{T} (\mathbb C)$ which connects the fixed points
$\ell_s$ and $r_t.$   This accounts for all the one dimensional orbits of
$\widetilde{T}(\mathbb C)$ in $Y.$  The torus $\widetilde{T}(\mathbb C)$ acts
on the orbit $E_{st}$ through the affine root $(\alpha, s+t).$ The orbit
$E_{st}$ is contained in $Y_{\gamma}$ iff $-v \le t-s \le v-1.$
\end{prop}

\subsection{Proof} If $s-t \ge 1$ then $E_{st}$  is the orbit of the point
$\ell_s x_{t-s}.$ It is contained in $Y_{\gamma}$ iff $s-t \le v.$   If $s-t
\le 0$ then $E_{st}$ is the orbit of the point $\ell_t y_{s-t}.$  It is
contained in $Y_{\gamma}$ iff $t-s \le v-1.$  \qed

\subsection{} As in \S \ref{sec-EqHoSL2} the co-root $\alpha^{\vee}$
determines isomorphisms $\mathbf S(\widetilde{\mathfrak t}) \cong \mathbb
C[x,t]$ and $\mathcal D(\widetilde{\mathfrak t}) \cong \mathbb
C[\partial_x,\partial_t].$   Fix $\gamma \in \mathfrak t(\mathfrak o)$ and let
$v = \val(\alpha'(\gamma)).$  Let $M_{ab}$ be the vector subspace of $\mathbb
C[\widetilde{W}] \otimes_{\mathbb C} \mathbf S(\widetilde{\mathfrak t})$ which
is spanned by elements
\[ (\ell_a - r_b) \otimes g_{ab}((a+b)x-2t) \]
where $g_{ab}$ are polynomials.  Then $M_{ab}$ is a $\mathcal
D(\widetilde{\mathfrak t})$ module.  Define $Q_v = \sum_{a,b} M_{ab}$ to be
the sum of those subspaces $M_{ab}$ such that  $-v \le b-a \le v-1.$ Then the
inclusion $\widetilde{W} \subset Y_{\gamma}$ induces a short exact sequence on
$T(\mathbb C)$-equivariant homology,
\begin{equation}\label{eqn-flagsequence}\begin{CD}
0 @>>> Q_v\left\{ \partial_t \right\} @>>> \mathbb C[\widetilde{W}]
\otimes_{\mathbb C}
\mathbf S(\mathfrak t) @>>> H_*^{T(\mathbb C)}(Y_{\gamma}) @>>>0
\end{CD}\end{equation}
where the {\it module of relations} $Q_v\left\{ \partial_t \right\}$ denotes
the submodule of $Q_v$ which is annihilated by the differential operator
$\partial_t.$

\begin{prop}\label{prop-flagrelations}
Fix $v \ge 1.$  The module $Q_v\left\{\partial_t\right\}$ of relations is
spanned by
\begin{equation}\label{eqn-moduleofrelations1}
  \sum_{d=1}^v  (1-\alpha^{\vee})^d \mathbb C[\widetilde{W}]
\otimes \mathbf S(\mathfrak t)\left\{ \partial_x^d \right\}\end{equation}
and
\begin{equation}\label{eqn-moduleofrelations2}
\sum_{d=1}^v (1-\alpha^{\vee})^{d-1} (1-w_{\alpha}) \mathbb C[\widetilde{W}]
\otimes \mathbf S(\mathfrak t) \left\{\partial_x^d \right\}.
\end{equation}
\end{prop}

\subsection{Proof}  This section is parallel to \S \ref{sec-combinatorial}.
Fix $a,b \in \mathbb Z$ and fix $d \ge 1.$  Define
\[ f_{a,b;d} \in \mathbb C[\widetilde{W}] \otimes_{\mathbb C} \mathbb C[x,t]
\]
to be the element
\[ \sum_{u=a}^{a+d-1} \sum_{v=b}^{b+d-1} (-1)^{u-v} \binom{d-1}{u-a}
\binom{d-1}{v- b} (\ell_u-r_v)\otimes \left( (u+v)x-2t\right)^{d-1}.\]
As in \S \ref{sec-combinatorial} we find that
\begin{equation} \label{eqn-fabd}
  f_{a,b;d} = (-1)^{a+b}(d-1)!(1-\alpha^{\vee})^{d-1}(\ell_{a} - r_{b})
\otimes x^{d-1} \in \ker(\partial_t).\end{equation}
Fix $v \ge 1.$  If $1\le d \le v$ then  $f_{m,m;d} \in Q_v$ and $f_{m+1,m;d}
\in Q_v.$ (These are the relations coming from the one dimensional orbits in
each of the two types of irreducible components of $Y^{\SL}_{\le d},$ cf.
Remark \ref{subsec-constellations}).  Define $K_v \subset \mathbb
C[\widetilde{W}] \otimes_{\mathbb C} \mathbb C[x,t]$ to be the span
\[ K_v = \sum_{m\in \mathbb Z}\sum_{d=1}^v \mathbb Cf_{m,m;d} + \mathbb
Cf_{m+1,m;d}.\]
  Then $K_v$ is a $\mathcal D(\mathfrak t) = \mathbb C[\partial_x]$ module, and
equation (\ref{eqn-fabd}) implies that $K_v \subset Q_v\left\{ \partial_t
\right\}.$ An argument similar to that of \S \ref{subsec-Pv} may be used to
show that in fact $K_v = Q_v \left\{ \partial_t \right\}.$  Since
\begin{align*}
f_{m,m;d} + f_{m+1,m,d} &= (d-1)!(1-\alpha^{\vee})^d \ell_m \otimes x^{d-1} \\
f_{m+1,m+1;d} + f_{m+1,m,d} &= (d-1)!(1-\alpha^{\vee})^d r_m \otimes x^{d-1}
\end{align*}
we see that $K_v$ is the sum
\[  \sum_{d=1}^v  (1-\alpha^{\vee})^d \mathbb C[\widetilde{W}]
\otimes \mathbf S(\mathfrak t)\left\{ \partial_x^d \right\}
+ \sum_{d=1}^v \mathbb C f_{0,0;d}.\]
The second set of relations may be replaced by the less efficient but more
symmetric collection (\ref{eqn-moduleofrelations2}) by noting that for any
$m \ge 0,$
\begin{align*}
\sum_{j=-m}^{m-1}(f_{j,j;d} + f_{j+1,j;d}) + f_{m,m;d} &=
(d-1)!(1-\alpha^{\vee})^{d-1}
(\ell_{-m} - r_m)\otimes x^{d-1} \\
\sum_{j=-m}^{m-1}(f_{j+1,j+1;d} + f_{j+1,j;d})-f_{m,m;d} &=
(d-1)!(1-\alpha^{\vee})^{d-1}
(r_{-m} - \ell_m)\otimes x^{d-1}.
\end{align*}
Since $K_v$ is a $\mathbb C[\partial_x]$ module, we conclude that for any $1
\le d \le v,$ \[ f_{0,0;d}\in (1-\alpha^{\vee})^{d-1}(1-w_{\alpha}) \mathbb
C[\widetilde{W}] \otimes_{\mathbb C} \mathbf S(\mathfrak t) \left\{
\partial_x^d \right\} \subset K_v.  \qed \]

\subsection{The semisimple rank one case}  As in \S \ref{subsec-H1} let
$H$ be a connected reductive complex algebraic group of semisimple rank one,
let $T \subset B \subset H$ be a Borel pair, let $I\subset H(F)$ be the
corresponding Iwahori subgroup, and let $Y^H = H(F)/I$ be the affine flag
manifold for $H.$  Let $\alpha, \alpha^{\vee}$ be the simple positive root and
co-root.  These determine a reflection $w_{\alpha} \in W$ in the (finite) Weyl
group for $H,$ a 1-dimensional sublattice $\langle \alpha^{\vee} \rangle$ in
the lattice $\Lambda^H$ of translations of $T,$ and a differential operator
$\partial_{\alpha} \in \mathcal D(\mathfrak t)$ of degree 1. Let
$\widetilde{W} = \Lambda \rtimes W$ be the extended affine Weyl group. The
mapping $\SL \to H$ determined by $\alpha$ induces an inclusion
$\phi_{\alpha}: Y^{\SL} \to Y^H.$  Let $\gamma \in \mathfrak t(\mathfrak o)$
be a regular element and set $v = \val \alpha'(\gamma).$  Then the affine flag
manifold for $H$ is the disjoint union
\[ Y^H = \coprod_{\ell \in \Lambda^{H}/\langle \alpha^{\vee} \rangle}
\ell \phi_{\alpha}(Y^{\SL})\]
and the affine Springer fiber for $\gamma$ is the disjoint union
\[
Y^H_{\gamma} = \coprod_{\ell \in \Lambda^{H}/\langle \alpha^{\vee} \rangle}
\ell \phi_{\alpha}(Y^{\SL}_{\le v})\]
which we may unambiguously denote by $Y^H_{\le v}.$ It follows, as in \S
\ref{subsec-proof-cohomology}, that the $T(\mathbb C)$ equivariant homology of
$Y^H_{\gamma}$ is again given by the following:

\begin{prop}\label{prop-ss1flag}  Suppose $H$ has semisimple rank one.  Fix $T
\subset B \subset H$ and fix $v \ge 1.$  Then the inclusion of the $T(\mathbb
C)$ fixed points in the affine Springer fiber $Y^H_{\le v}$ induces a short
exact sequence
\begin{equation}\label{eqn-ss1flagsequence} \begin{CD}
0 @>>> Q_v\left\{ \partial_t \right\} @>>> \mathbb C[\widetilde{W}]
\otimes_{\mathbb C}
\mathbf S(\mathfrak t) @>>> H_*^{T(\mathbb C)}(Y^H_{\le v}) @>>>0
\end{CD}\end{equation}
where the module of relations $Q_v\left\{ \partial_t \right\}$ is the span of
\parens{{\rm\ref{eqn-moduleofrelations1}}} and
\parens{{\rm\ref{eqn-moduleofrelations2}}}.
\qed
\end{prop}

\section{Homology of affine Springer fibers in the affine flag manifold}
\label{sec-flag}
The results in this section are parallel to \S \ref{sec-equivhomology} and \S
\ref{sec-endoscopic} so they will be presented without detailed proofs.
\subsection{}
As in \S \ref{sec-affinespringer} let $F = \mathbb C((\epsilon)),$ let $G$ be
a connected reductive complex algebraic group, let $T \subset B \subset G$ be
a Borel pair and let $I \subset G(F)$ be the corresponding Iwahori subgroup.
The affine flag manifold is $Y = Y^G = G(F)/I.$  Let $\Lambda \subset T(F)$
denote the lattice of translations.  Identify each co-root $\alpha^{\vee} \in
\Phi^{\vee}(G,T)$ with its image $\alpha^{\vee}(\epsilon) \in \Lambda.$  Set
$\widetilde{W} = \Lambda \rtimes W$ where $W = W(G(\mathbb C),T(\mathbb C))$
is the Weyl group. Each cell in the Bruhat decomposition
\[ Y = \coprod_{w\in \widetilde{W}} IwI/I \]
contains a unique $T(\mathbb C)$ fixed point.  Hence the choice of basepoint
determines a one to one correspondence
\[ Y^{T(\mathbb C)} \cong \widetilde{W} \]
between the set of fixed points and the extended affine Weyl group, which is
compatible with the action of the lattice $\Lambda$ of translations.

For each root $\alpha \in \Phi^+$ let $w_{\alpha}\in W$ be its corresponding
reflection, let $W_{\alpha} = \left\{ 1, w_{\alpha} \right\},$ and let
$H^{\alpha}$ be the connected reductive group of semisimple rank one which
contains $T$ and the root subgroup $U_{\alpha} \subset G.$ Let $Y^{\alpha}$ be
the affine flag manifold for $H^{\alpha}.$ If $u \in Y$ is a $T(\mathbb C)$
fixed point then it determines an isomorphism
\[ \phi_u: Y^{\alpha} \to H^{\alpha}\cdot u \subset Y \]
between $Y^{\alpha}$ and the orbit $H^{\alpha} \cdot u.$  It restricts to an
isomorphism $Y^{\alpha}_{\le v} \cong Y_{\gamma} \cap H^{\alpha}\cdot u$ of
affine Springer fibers, for any regular element $\gamma \in \mathfrak
t(\mathfrak o)$ with $v = \val \alpha'(\gamma).$  (The action of $w_{\alpha}$
on $Y^{T(\mathbb C)}$ preserves the fixed points in the orbit $H^{\alpha}
\cdot u$ but it interchanges those of ``type $\ell$'' and ``type $r$''.)
Define
\[ Z^{\alpha} = \coprod_{u \in W_{\alpha} \backslash W} H^{\alpha} \cdot u =
\coprod_{u \in W_{\alpha} \backslash W} \phi_u(Y^{\alpha}).\]
The proof of the following lemma is parallel to that of Lemma
\ref{lem-1dimensional}.

\begin{lem}  Let $\gamma \in \mathfrak t(\mathfrak o)$ be a regular element.
Then the union of the 0 and 1 dimensional orbits in the affine Springer fiber
is
\[ (Y_{\gamma})_1 = \bigcup_{\alpha \in \Phi^+} Z^{\alpha}_{\gamma} \]
where $Z^{\alpha}_{\gamma}$ is the intersection
\[ Z^{\alpha}_{\gamma} = Y_{\gamma} \cap Z^{\alpha} = \coprod_{u\in W_{\alpha}
\backslash W} \phi_u(Y^{\alpha}_{\le v}). \]
If $\alpha \ne \beta$ then $Z^{\alpha}_{\gamma} \cap Z^{\beta}_{\gamma} =
\widetilde{W}.$  \qed
\end{lem}

Each $\alpha \in \Phi^+$ corresponds to a degree one differential operator
$\partial_{\alpha} \in \mathcal D(\mathfrak t),$ to an element $\alpha^{\vee}
\in \Lambda,$ and to a reflection $w_{\alpha}\in W.$  For $\alpha \in \Phi^{+}
$ define the $\mathcal D(\mathfrak t)$ submodule $M_{\alpha,\gamma}$ to be the
sum
\begin{equation*}
\sum_{d=1}^{\val(\alpha'(\gamma))} (1-\alpha^{\vee})^d \mathbb
C[\widetilde{W}] \otimes_{\mathbb C} \mathbf S(\mathfrak t) \left\{
\partial_{\alpha}^d \right\}
+ \sum_{d=1}^{\val(\alpha'(\gamma))} (1-\alpha^{\vee})^{d-1}(1-w_{\alpha})
\mathbb C[\widetilde{W}] \otimes_{\mathbb C} \mathbf S(\mathfrak t) \left\{
\partial_{\alpha}^d\right\}.
\end{equation*}
Using Proposition \ref{prop-ss1flag}, the argument of \S
\ref{subsec-proof-cohomology} gives the following.

\begin{thm}\label{thm-affineflaghomology}
Let $\gamma \in \mathfrak t(\mathfrak o)$ be a regular element. Suppose the
homology $H_*(Y_{\gamma}; \mathbb C)$ is pure.  Then the inclusion
$\widetilde{W} \subset Y_{\gamma}$ induces an exact sequence of $\mathcal
D(\mathfrak t)$ modules,
\begin{equation*} \begin{CD}
0 @>>> \sum_{\alpha\in \Phi^+} M_{\alpha,\gamma} @>>> \mathbb C[\widetilde{W}]
\otimes_{\mathbb C} \mathbf S(\mathfrak t) @>>> H_*^{T(\mathbb C)}(Y_{\gamma})
@>>> 0.
\end{CD}  \end{equation*}
\end{thm}
The group $(\widetilde{W} \rtimes \mathbf{Aut})_{\gamma}$ of
(\ref{eqn-Wgamma}) acts on this equivariant homology group and it restricts to
an action on the ordinary homology,
\[ H_*(Y_{\gamma}) = H^{T(\mathbb C)}_*(Y_{\gamma})\left\{ \mathcal I
\right\}, \]
that is, the subgroup of the equivariant homology which is annihilated by the
augmentation ideal $\mathcal I \subset \mathcal D.$

\subsection{Springer action}\label{subsec-Springeraction}
   Consider the (regular$\otimes$trivial) action of $\widetilde{W}$ from the
right on 
\[\mathbb C[\widetilde{W}] \otimes \mathbf S(\mathfrak t).\]  
This
action preserves each relation $M_{\alpha,\gamma}$ (for $\alpha \in \Phi^+$).
Assume the homology $H_*(Y_{\gamma}; \mathbb C)$ is pure.  Then the right
action of $\mathbb C[\widetilde{W}]$  passes to an action on $H_*^{T(\mathbb
C)}(Y_{\gamma})$  which is compatible with the $\mathcal D(\mathfrak t)$
module structure, and which commutes with the $(\widetilde{W} \rtimes
\mathbf{Aut})_{\gamma}$ action.  So it restricts to an action (from the right)
of $\widetilde{W}$ on the ordinary homology $H_*(Y_{\gamma}).$  It can be
shown that this action coincides with the Springer action defined by 
Lusztig \cite{lusztig} and (using a statement from \cite{KL0} whose proof
does not appear in the literature) by Sage \cite{Sage2}. So  Theorem 
\ref{thm-affineflaghomology} gives a ``formula'' for the Springer action.

\subsection{Endoscopic groups}  There are parallels to \S \ref{sec-endoscopic}
and \S \ref{sec-quotient} for the affine flag manifold also.  As in \S
\ref{sec-endoscopic} suppose $G$ is adjoint, and let $(H,s) $ be endoscopic
data for $G.$  Fix compatible Borel pairs $T \subset B \subset G$ and $T_H
\subset B_H \subset H$  with resulting identification $T_H \cong T.$
The group $H$ has root system $\Phi(H,T_H)\subset \Phi(G,T)$ and (extended)
affine Weyl group $\widetilde{W}^H = \Lambda \rtimes W^H\subset
\widetilde{W}.$  Let $\gamma \in \mathfrak t(\mathfrak o)$ be a regular
element with corresponding Springer fibers $Y_{\gamma}$ and $Y^H_{\gamma_H}.$
Suppose the homology groups $H_*(Y_{\gamma}; \mathbb C)$ and
$H_*(Y^H_{\gamma_H};\bC)$ are pure.   Let $\Delta$ be the homological transfer
factor of (\ref{eqn-transfer}), let $r = \deg(\Delta)$, and (as
before) let $\eta$ be the sign
character that gives the action of  $\widetilde{\mathcal
W}^{G,H}_{\gamma}$ on $\Delta$.

We are going to compare  the induced module
$H_*^{T(\mathbb C)}(Y^H_{\gamma_H})[-2r]\otimes_{\bC[\widetilde
W^H]}\bC[\widetilde W]$ to
the module $ H_*^{T(\mathbb C)} (Y_{\gamma})$. One sees immediately
that this induced module is
equal to the quotient of $\mathbb C[\widetilde{W}] \otimes_{\mathbb C} \mathbf
S(\mathfrak t) $ by $\sum_{\alpha\in \Phi^+(H)} M_{\alpha,\gamma}$.
Note that the surjection
\[
1 \otimes \Delta:\mathbb C[\widetilde{W}] \otimes_{\mathbb C} \mathbf
S(\mathfrak t) \to \mathbb C[\widetilde{W}] \otimes_{\mathbb C} \mathbf
S(\mathfrak t)
\]
kills $M_{\alpha,\gamma}$ whenever $\alpha \in \Phi^+(G) - \Phi^+(H).$

  Let $J \subset\mathbb C[\Lambda]$  be the
multiplicative subset generated by the elements $1 - \alpha^{\vee}$ for
$\alpha^{\vee} \in \Phi^{\vee}(G) - \Phi^{\vee}(H).$ We use the
action of $\Lambda$ on the
Springer fibers to regard their homology groups as modules over the
group algebra $\mathbb
C[\Lambda]$, allowing us to localize the homology groups using $J$.
Then we have the following
result.

\begin{thm} \label{thm-FLforflags}  Suppose $H_*(Y_{\gamma})$ and
$H_*(Y^H_{\gamma_H})$ are pure.  Then the mapping $1\otimes\Delta$
induces a homomorphism of $\mathcal D(\mathfrak t)$ modules,
\[ H_*^{T(\mathbb C)} (Y_{\gamma}) \to H_*^{T(\mathbb C)}(Y^H_{\gamma_H})[-2r]
\otimes_{\bC[\widetilde W^H]}\bC[\widetilde W]
\]
and a homomorphism on ordinary homology,
\[H_*(Y_{\gamma}; \mathbb C) \to H_*(Y^H_{\gamma_H}; \mathbb C)[-2r]
\otimes_{\bC[\widetilde W^H]}\bC[\widetilde W].  \]
Both homomorphisms transform under $\eta$ with respect to the left action of
$\widetilde{\mathcal W}^{G,H}_{\gamma}$ and are equivariant with
respect to the right
action of $\widetilde{W},$ and both homomorphisms become isomorphisms
\begin{gather*}
H_*^{T(\mathbb C)} (Y_{\gamma})_J \cong H_*^{T(\mathbb C)}
(Y^H_{\gamma_H})_J[-2r]\otimes_{\bC[\widetilde W^H]}\bC[\widetilde W] \\
  H_*(Y_{\gamma}; \mathbb C)_J \cong H_*(Y^H_{\gamma_H}; \mathbb C)_J[-2r]
\otimes_{\bC[\widetilde W^H]}\bC[\widetilde W]
\end{gather*}
after localizing with respect to $J.$\qed \end{thm}

\section{An orbital integral}\label{sec-orbitalintegral}
\subsection{Set-up}\label{1}
Let $k$ be a finite field and $\bar k$ an algebraic closure of~$k$.
Let $F=k((\varepsilon))$ and
$L=\bar k((\varepsilon))$. Let $\mathfrak o$ denote the valuation
ring $k[[\varepsilon]]$ of~$F$.
We write
$\sigma$ for the Frobenius automorphism of both
$\bar k/k$ and
$L/F$, and we write $\langle\s\rangle$ for the infinite cyclic group
generated by~$\s$. We also
choose an algebraic closure
$\bar L$ of~$L$ and write $I$ for
$\Gal(\bar L/L)$. We write $\Fbar$ for the algebraic closure of~$F$
in~$\Lbar$ and $\G$ for the
Galois group $\Gal(\Fbar/F)$.

Let $G$ be a connected reductive group over~$F$, and assume that $G$
splits over~$L$.
Let $K$ be a parahoric subgroup of~$G(F)$, let $K_L$ be the
corresponding $\s$-stable parahoric
subgroup of~$G(L)$, and  write
$\mathfrak k$ and $\mathfrak k_L$ for the corresponding
parahoric subalgebras (of $\mathfrak g(F)$ and
$\mathfrak g(L)$ respectively). Let $X$ be the
$k$-ind-scheme  associated to $G(L)/K_L$. Then $X(\bar k)=G(L)/K_L$
and $X(k)=G(F)/K$.

Let $T$ be a maximal $F$-torus in~$G$. Let $S$ denote the maximal
unramified subtorus of~$T$. Thus
the cocharacter group $X_*(S)$ can be identified with $X_*(T)^I$.

Let $u$ be a regular element in $\Lie(T)(F)$ and assume that $u$ is
\emph{integral}, in the sense
that $\alpha(u)$ lies in the valuation ring of $\bar F$ for every root $\alpha$
of~$T$ in~$G$.
We
write $X_u$ for the affine Springer fiber
$\{x \in X: x^{-1}ux \in \mathfrak k_L\}$ studied by
Kazhdan-Lusztig in \cite{KL}. (We usually write $x^{-1}ux$ rather than
$\Ad(x)^{-1}u$.)

For $\mu \in X_*(S)$ we write $\varepsilon^\mu$ for $\mu(\varepsilon)
\in S(L)$. The map $\mu
\mapsto \varepsilon^\mu$  is $\langle\s\rangle$-equivariant and
identifies $X_*(S)$ with a subgroup
of~$S(L)$ (and of~$T(L)$) that we will denote by~$\L$.  The group
$T(L)$ centralizes $u$ and
therefore acts by left translations on~$X_u$. The quotient
$\L\backslash X_u$ of $X_u$ by the
subgroup $\L$ of~$T(L)$ is a non-empty projective scheme of finite
type over~$k$
  (see \cite{KL}).

\subsection{Normalized Haar measure  on $T(F)$}\label{2}
Recall from \cite[7.6]{Kottwitz} the exact sequence
\begin{equation}\label{2.1}
1 \to T(L)_1 \to T(L) \xrightarrow{w_T} X_*(T)_I \to 0.
\end{equation}
Taking invariants under $\langle\s\rangle$, we get another exact sequence
\begin{equation}\label{2.2}
1 \to T(F)_1 \to T(F) \to (X_*(T)_I)^{\langle\s\rangle} \to 0,
\end{equation}
where $T(F)_1:=T(F)\cap T(L)_1$. Surjectivity at the right end of
this last exact sequence is
proved in \cite[7.6]{Kottwitz}.

Let $dt$ be the Haar measure on~$T(F)$ that gives $T(F)_1$
measure~$1$. The group
$\L^{\langle\s\rangle}$ of~$T(F)$ is discrete and cocompact, so the
volume of the quotient
$\L^{\langle\s\rangle}\backslash T(F)$ is finite. We need to compute
this volume (with
respect to $dt$).

\subsection{Volume computation}\label{3}
We claim that
\begin{equation}\label{3.1}
\vl_{dt}(\L^{\langle\s\rangle}\backslash T(F))=\frac{|\cok[X_*(S)_\G
\to X_*(T)_\G]|}
{|\ker[X_*(S)_\G \to X_*(T)_\G]|}.
\end{equation}
To prove this claim we begin by noting that the canonical map
\[X_*(S)= X_*(T)^I \to X_*(T)_I\]
has finite kernel and cokernel. Thus this map is in fact injective,
since $X_*(S)$ is
torsion-free. We let $M$ denote its cokernel, so that we get a short
exact sequence
\[
0 \to X_*(S) \to X_*(T)_I \to M \to 0
\]
of $\langle\s\rangle$-modules, and an associated long exact sequence
of group cohomology
\begin{equation}\label{3.2}
\begin{split}
0 \to &X_*(S)^{\langle\s\rangle} \to (X_*(T)_I)^{\langle\s\rangle}
\to M^{\langle\s\rangle}\\
  \to  &X_*(S)_{\G}\ \  \to \ \ X_*(T)_\G \quad\  \to
M_{\langle\s\rangle} \to 0,
\end{split}
\end{equation}
in which we have identified $X_*(S)_{\langle\s\rangle}$ with $X_*(S)_{\G}$  and
$(X_*(T)_I)_{\langle\s\rangle}$ with $X_*(T)_{\G}$.

Since $X_*(S)^{\langle\s\rangle} \to (X_*(T)_I)^{\langle\s\rangle}$
is injective, the subgroups
$X_*(S)^{\langle\s\rangle}$ and $T(F)_1$ of~$T(F)$ have trivial
intersection, and therefore
\begin{equation*}
\vl_{dt}(\L^{\langle\s\rangle}\backslash
T(F))=\left|\frac{(X_*(T)_I)^{\langle\s\rangle}}{X_*(S)^{\langle\s\rangle}}
\right|.
\end{equation*}
It follows from \eqref{3.2} that
\begin{equation*}
\left |\frac{(X_*(T)_I)^{\langle\s\rangle}}{X_*(S)^{\langle\s\rangle}}\right
|=\frac{|M^{\langle\s\rangle}|}{|\ker[X_*(S)_{\G}  \to X_*(T)_\G]|}.
\end{equation*}
Since all groups in the exact sequence
\[
0 \to M^{\langle\s\rangle} \to M \xrightarrow{\s -1} M \to
M_{\langle\s\rangle} \to 0
\]
are finite, we have
\begin{equation*}
|M^{\langle\s\rangle}|=|M_{\langle\s\rangle}|.
\end{equation*}
Finally, it follows from  \eqref{3.2} that
\begin{equation*}
|M_{\langle\s\rangle}|=\cok[X_*(S)_{\G}  \to X_*(T)_\G].
\end{equation*}
Combining the last four equations, we get the claim.

\subsection{Normalized orbital integrals}\label{4}
For any compactly supported locally constant function $f$ on
$\mathfrak g(F)$ we put
\[
O_u(f)=\int_{T(F)\backslash G(F)} f(g^{-1}ug) \, dt\backslash dg,
\]
where $dt$ is the normalized Haar measure on~$T(F)$ defined above,
and $dg$ is the Haar measure
on~$G(F)$ that gives our chosen parahoric subgroup $K$ measure $1$.
We then have
\begin{equation}\label{4.1}
O_u(f)=\vl_{dt}(\L^{\langle\s\rangle}\backslash T(F))^{-1}
\int_{\L^{\langle\s\rangle}\backslash
G(F)} f(g^{-1}ug) \, dg.
\end{equation}
  In the special case that our function $f$ is $\mathbf 1_{\mathfrak k}$, the
characteristic
function of $\mathfrak k$, we have
\begin{equation}\label{4.2}
\begin{split}
  \int_{\L^{\langle\s\rangle}\backslash
G(F)} \mathbf 1_{\mathfrak k}(g^{-1}ug) \, dg &=|\{x \in
\L^{\langle\s\rangle}\backslash G(F)/K :
x^{-1}ux \in \mathfrak k \}| \\
&=|\L^{\langle\s\rangle}\backslash(X_u(k))|.
\end{split}
\end{equation}
Combining \eqref{3.1}, \eqref{4.1}, \eqref{4.2}, we find that
\begin{equation}\label{4.3}
O_u(\mathbf 1_{\mathfrak k})=\frac{|\ker[X_*(S)_\G \to X_*(T)_\G]|}
{|\cok[X_*(S)_\G \to X_*(T)_\G]|}\cdot
|\L^{\langle\s\rangle}\backslash(X_u(k))|.
\end{equation}

\subsection{$\kappa$-orbital integrals of $\mathbf 1_{\mathfrak k}$}\label{5}
Consider an element $u' \in \mathfrak g(F)$ that is stably conjugate
to~$u$ (that is, conjugate
under $G(\Fbar)$). Since the group $H^1(L,T)$ vanishes (see
\cite[X.7]{serre68}), there exists $ g
\in G(L)$ such that $u' = gug^{-1}$. Apply $\s$ to this equality to
see that $g^{-1}\s(g) \in
T(L)$. The class of $g^{-1}\s(g)$ in~$B(T)$ (in other words, the
$\s$-conjugacy class of
$g^{-1}\s(g)$ in~$T(L)$) is independent of the choice of~$g$. Here,
as usual, for any linear
algebraic group $G$ over $F$ we write $B(G)$ for the set of
$\sigma$-conjugacy classes in
$G(L)$.

Turning this around, given $t \in T(L)$ whose class in~$B(T)$ lies in
\[
\mathfrak D(T/F):=\ker[B(T) \to B(G)],
\]
we choose $g \in G(L)$ such that $t = g^{-1}\s(g)$ and put
$u(t)=gug^{-1}$, an element
of~$\mathfrak g(F)$ that is stably conjugate to~$u$ and whose
$G(F)$-conjugacy class depends only
on the class of~$t$ in~$B(T)$. The construction $t \mapsto u(t)$ sets
up a bijection from
$\mathfrak D(T/F)$ to the set of $G(F)$-conjugacy classes in the
stable conjugacy class of~$u$.

Recall from \cite[7.6]{Kottwitz} that the canonical homomorphism
\[
w_T:T(L) \to X_*(T)_I
\]
induces an isomorphism (take coinvariants under $\langle\s\rangle$)
\[
B(T)\simeq X_*(T)_\G.
\]

At this point we fix a prime number $l$ that is non-zero in the field $k$.
Put
\[\hat T:=\Hom(X_*(T),\qlbarstar).
\]
  Let $\kappa \in \Hom(X_*(T)_\G,\qlbarstar)=\hat T^\G$. For
$t \in T(L)$ we write $\langle t, \kappa \rangle \in \qlbarstar$ for
the value of the
character~$\kappa$ on the element of $X_*(T)_\G$ obtained as the
image of~$t$ under
\[
T(L) \twoheadrightarrow B(T) \simeq X_*(T)_\G.
\]

For $\Qbar_l$-valued locally constant compactly supported functions
$f$ on $\mathfrak g(F)$ we
define the
$\kappa$-orbital integral
$O^\kappa_u(f)$ by
\[
O^\kappa_u(f)=\sum_{t \in \mathfrak D(T/F)} \langle t,\kappa \rangle
\cdot O_{u(t)}(f).
\]
Note that if $\kappa$ lies in the subgroup $Z(\hat G)^\G$ of $\hat
T^\G$ (where we form
the Langlands dual group $\hat G$ using $\Qbar_l$, and where $Z(\hat
G)$ denotes the center of
$\hat G$), then $O^\kappa_u$ is a stable orbital integral.

   In the special case that
$f$ is
$\mathbf 1_{\mathfrak k}$, we find that
\begin{equation}\label{5.1}
O^\kappa_u(\mathbf 1_{\mathfrak k})
=\frac{|\ker[B(S)\to B(T)]|}{|\cok[B(S)\to B(T)]|} \cdot \sum_{t \in
\mathfrak D(T/F)}
\langle t,\kappa \rangle \cdot  |\L^{\langle\s\rangle}\backslash X_u^{t\s}|,
\end{equation}
where we have written $X_u^{t\s}$ for the fixed point set of $t\s$ on
$X_u$, which appears in our
formula since $x \mapsto gx$ induces a bijection from $X_u^{t\s}$ to
$X_{u(t)}(k)$, where $g \in
G(L)$ is chosen so that $t=g^{-1}\s(g)$, as before.

\subsection{Local systems $\mathcal L_\eta$ on $\L\backslash X_u$}
\label{subsec-localsystem}
Any finite dimensional (continuous) $l$-adic representation of $\L
\rtimes \Gal(\bar k/k)$
(with the topology on $\L$ given by all subgroups of finite index)
gives
rise to a local system on $\L\backslash X_u$ over~$k$ (which becomes
constant when pulled back to
$X_u$ over $\bar k$). Put $\hat S:=\Hom(X_*(S),\qlbarstar)$. Let
\[
\eta \in \Hom(X_*(S)_\G,\qlbarstar)=\hat S^\G=\hat
S^{\langle\s\rangle},
\]
and assume that $\eta$ has finite order.
Using the isomorphism $\L \simeq X_*(S)$, we also view $\eta$ as a
character on~$\L$, and since
$\eta$ is fixed by~$\s$, we may extend $\eta$ to a character on $\L \rtimes
\Gal(\bar k/k)$ by
making it trivial on $\Gal(\bar k/k)$. We use the
resulting $1$-dimensional representation of
$\L \rtimes \Gal(\bar k/k)$ to form a rank $1$ local system $\mathcal
L_\eta$ on $\L\backslash
X_u$ over~$k$.

Since the action of~$T(L)$ on~$X_u$ commutes with that of~$\L$, it
induces an action (over $\bar
k$) of~$T(L)$ on $\L\backslash X_u$ and $\mathcal L_\eta$, and
therefore $T(L)$ acts on the
cohomology groups $H^i(\L\backslash X_u,\mathcal L_\eta)$. It follows
from the definitions that
the subgroup $\L$ of~$T(L)$ acts on these cohomology groups through
the character $\eta^{-1}$.
Moreover the ``identity component'' $T(L)_1$ of~$T(L)$ acts
trivially on them. Then since $\L \cdot T(L)_1$ has finite index
in~$T(L)$ and acts
by a character on cohomology, each cohomology group is a semisimple
$T(L)$-module. Therefore
(using \eqref{2.1}) we can decompose the cohomology groups as
\[
  H^i(\L\backslash X_u,\mathcal L_\eta)=\bigoplus_\kappa
H^i(\L\backslash X_u,\mathcal L_\eta)_\kappa
\]
where $\kappa$ runs through the finite set of characters $\kappa \in
\Hom(X_*(T)_I,\qlbarstar)=\hat
T^I$ whose image under $\hat T^I \to \hat S$ is $\eta^{-1}$, and
where $H^i(\L\backslash
X_u,\mathcal L_\eta)_\kappa$ denotes the $\kappa$-isotypic subspace
of $H^i(\L\backslash
X_u,\mathcal L_\eta)$ (viewing $\kappa$ as a character on~$T(L)$, as before).

\subsection{Main result}\label{subsec-mainresult}
Now we can state the main result of this section, which gives a
cohomological interpretation of the
$\kappa$-orbital integrals of~$\mathbf 1_{\mathfrak k}$.
\begin{thm}\label{orb-thm}
Let $\kappa$ be an element of finite order in $\hat T^\G$.
  Then the $\kappa$-orbital
integral
$O^\kappa_u(\mathbf 1_{\mathfrak k})$ is given by
\[
O^\kappa_u(\mathbf 1_{\mathfrak k})=
\textup{Tr}(\s^{-1};H^\bullet(\L\backslash X_u,\mathcal L_\eta)_\kappa),
\]
where $\eta$ is the
image of $\kappa^{-1}$ under $\hat T^\G \to \hat S^\G$.
\end{thm}
Here we have written $\text{Tr}(\s^{-1};H^\bullet(\L\backslash
X_u,\mathcal L_\eta)_\kappa)$ as an
abbreviation for
\[
\sum^{2\dim X_u}_{i=0} (-1)^i \text{Tr}(\s^{-1};H^i(\L\backslash
X_u,\mathcal L_\eta)_\kappa).
\]
\begin{proof}
The first step is to apply the Grothendieck-Lefschetz trace formula.
We claim that for any $ t \in
T(L)$ we have
\begin{equation}\label{7.1}
\text{Tr}\bigl((t\s)^{-1};H^\bullet(\L\backslash X_u,\mathcal
L_\eta)\bigr)=\sum_{\bar x \in
(\L\backslash X_u)^{t\s}} \langle \l_x, \eta \rangle,
\end{equation}
with notation as follows. We are writing $(\L\backslash
X_u)^{t\s}$ for the fixed points of $t\s$ on $\L\backslash X_u$, and for such a
fixed point $\bar x$ we choose a representative $x \in X_u$ and
define $\l_x \in \L$ by the
equality $t\s x=\l_x x$. The image of $\l_x$ in
$\L_{\langle\s\rangle}=X_*(S)_\G$ is independent
of the choice of representative~$x$. The equality \eqref{7.1} follows from the
Grothendieck-Lefschetz trace formula applied to the Frobenius map for a twisted
$k$-form of $\L\backslash X_u$, $\mathcal L_\eta$ in which the
twisted action of $\s$ is given by
$t\s$.

The formula \eqref{7.1} can be rewritten as
\begin{equation}\label{7.2}
\sum_{\kappa'} \text{Tr}\bigl((t\s)^{-1};H^\bullet(\L\backslash X_u,\mathcal
L_\eta)_{\kappa'}\bigr)=\sum_{\l \in X_*(S)_\G} \langle \l, \eta
\rangle^{-1} \cdot
|\L^{\langle\s\rangle}\backslash (X_u^{\varepsilon^\l t\s})|,
\end{equation}
where the sum on the left is taken over $\kappa' \in \hat T^\G$ whose
image under $\hat T^\G \to
\hat S^\G$ is $\eta^{-1}$. Indeed, since $(t\s)^{-1}$ maps $H^i(\L
\backslash X_u,\mathcal
L_\eta)_{\kappa'}$ to  $H^i(\L \backslash X_u,\mathcal
L_\eta)_{\s^{-1}(\kappa')}$,
only those $\kappa' \in \hat T^I$ that are fixed by~$\sigma$
contribute to the trace. As for the
sum on the right, we simply collected like terms in the sum on the
right in \eqref{7.1}.

Multiply both sides of \eqref{7.2} by $\langle t,\kappa \rangle$, sum over
\[ t \in \cok [B(S) \to B(T)]=\cok[X_*(S)_\G \to X_*(T)_\G],
\]
and divide by the order $|\cok|$ of this cokernel. We obtain
\[
\text{Tr}\bigl(\s^{-1};H^\bullet(\L\backslash X_u,\mathcal
L_\eta)_{\kappa}\bigr)=\frac{|\ker|}{|\cok|} \cdot \sum_{t \in B(T)}
\langle t, \kappa \rangle \cdot
|\L^{\langle\s\rangle}\backslash (X_u^{t\s})|.
\]
Comparing this equation with equation \eqref{5.1}, we see that in
order to finish the proof of the
theorem, it suffices to show that if $X_u^{t\s}$ is non-empty, then
$t$ lies in $\ker[B(T) \to
B(G)]$. Indeed, suppose that $t\s(x)=x$ for some $x \in G(L)/K_L$.
Then $x^{-1}t\s(x) \in K_L$. By
Lang's theorem (applied to the finite dimensional quotients of the
connected proalgebraic group
$K_L$ over $k$) and the completeness of~$L$, every element in~$K_L$
is trivial in~$B(G)$.
Therefore $t$ is trivial in~$B(G)$.
\end{proof}

\subsection{A variant of the main result}\label{8}
The theorem above can be restated in a slightly different way, since
the classical
relationship between the cohomology of $X_u$ and that of
$\Lambda\backslash X_u$
presumably has an
$l$-adic version
\begin{equation*}
H^i(\L\backslash X_u,\mathcal
L_\eta)_\kappa=Ext^i_{X_*(T)_I}((\Qbar_l)_\kappa,R\G(X_u,\Qbar_l)).
\end{equation*}
Here $(\Qbar_l)_\kappa$ denotes the $1$-dimensional $X_*(T)_I$-module
obtained by letting
$X_*(T)_I$ act on $\Qbar_l$ via the surjection $X_*(T)_I \to X_*(T)_\G$ and the
character $\kappa$ of $X_*(T)_\G$, and the complex $R\G(X_u,\Qbar_l)$ is being
regarded as an object in a suitable
  derived category of $X_*(T)_I$-modules. Note that $\sigma$ acts on
$X_*(T)_I$ as well as $R\G(X_u,\Qbar_l)$; it then acts on
\[
Ext^i_{X_*(T)_I}((\Qbar_l)_\kappa,R\G(X_u,\Qbar_l))
\]
as well since $\kappa$ is fixed by
$\sigma$. Therefore the theorem above presumably yields the equality
\begin{equation}\label{variant.main}
O^\kappa_u(\mathbf 1_{\mathfrak k})=
\text{Tr}(\s^{-1};RHom^\bullet_{X_*(T)_I}((\Qbar_l)_\kappa,R\G(X_u,\Qbar_l))).
\end{equation}

\subsection{A geometric reformulation of the fundamental lemma}\label{9}
We now assume that $G$ comes from a (necessarily quasi-split)
connected reductive group over $k$,
which we still denote by~$G$. Thus $G$ is unramified over $F$ and the
parahoric subgroup
$K:=G(\mathfrak o)$ is a hyperspecial maximal compact subgroup of
$G(F)$; we denote the
corresponding parahoric subalgebra by $\mathfrak k:=\mathfrak g(\mathfrak o)$.

As above we form all Langlands dual groups using $\Qbar_l$. Let $H$
be an endoscopic group for $G$
and let $s$ be the usual element in $Z(\hat H)^\G$, where $Z(\hat H)$
denotes the center of $\hat
H$. We write $\mathfrak h$ for the Lie algebra of $H$. We assume that
$H$ is also unramified, so
that it comes from a (quasi-split) group over
$k$, which we still denote by $H$. Thus we also have $K_H:=H(\mathfrak o)$ and
$\mathfrak k_H:=\mathfrak h(\mathfrak o)$.

We use a regular nilpotent element in $\mathfrak g(\mathfrak o)$
whose image in $\mathfrak
g(k)$ is also regular nilpotent in order to form Kostant's section,
obtaining as in
\cite{kottwitz98c} transfer factors $\Delta(u_H,u)$ which for $u$
lying in Kostant's section take
the simple value
\[
D_G(u)D_H(u_H)^{-1},
\]
(in other words, $\Delta_{IV}(u_H,u)$ in the notation of
Langlands-Shelstad \cite{ls}).
Here $D_G$ denotes the usual discriminant function on the Lie algebra.

Let $u_H$ be an integral $G$-regular semisimple element
of $\mathfrak h(F)$ and let $u$ be an image of $u_H$ (in the sense of
Langlands-Shelstad \cite{ls})
that is $G(F)$-conjugate to an element in Kostant's section. (This
uniquely determines the
$G(F)$-conjugacy class of $u$.) We then have affine Springer fibers
$X_u$, $X^H_{u_H}$ for $u$,
$u_H$ respectively. Using the dimension formula for affine Springer
fibers (conjectured by
Kazhdan-Lusztig \cite{KL} and proved by Bezrukavnikov \cite{Bez}) we
can rewrite the transfer factor
as
\[
\Delta(u_H,u)=q^{-r}
\]
where $r:=\dim X_u - \dim X^H_{u_H}$ and $q$ denotes the cardinality of~$k$.

We let $T$ denote the centralizer of $u$ in~$G$, and we identify $T$
with the centralizer $T_H$ of
$u_H$ in~$H$ via the unique admissible isomorphism whose differential
maps $u_H$ to~$u$. Via the
canonical injection $Z(\hat H) \hookrightarrow \hat T$ the element $s
\in Z(\hat H)^\G$ determines
an element $\kappa \in \hat T^\G$.

The conjectural fundamental lemma predicts the equality
\[
O^\kappa_{u_H}(\mathbf 1_{\mathfrak k_H})=q^{-r} \cdot
O^\kappa_u(\mathbf 1_{\mathfrak k}),
\]
which in view of \eqref{variant.main} can be rewritten as the equality of
\[
q^r \cdot
\text{Tr}(\s^{-1};RHom^\bullet_{X_*(T)_I}((\Qbar_l)_\kappa,R\G(X^H_{u_
H},\Qbar_l)))
\]
and
\[
\text{Tr}(\s^{-1};RHom^\bullet_{X_*(T)_I}((\Qbar_l)_\kappa,R\G(X_u,\Qbar_l)))
\]

We conjecture the stronger statement that
\[
RHom^\bullet_{X_*(T)_I}((\Qbar_l)_\kappa,R\G(X^H_{u_H},\Qbar_l[-2r](-r)))
\]
is isomorphic to
\[
RHom^\bullet_{X_*(T)_I}((\Qbar_l)_\kappa,R\G(X_u,\Qbar_l)).
\]

\subsection{Transfer factors for unramified tori}\label{10}
  Earlier in this paper, when
we calculated the homology of Springer fibers, the elements of
$\mathfrak g(F)$ we considered did
not necessarily lie in Kostant's section. Therefore we need a further
discussion of transfer
factors.

We continue the discussion in the subsection above, making only the
following changes. We no
longer require that the elements in $\mathfrak g(F)$ that we consider
be $G(F)$-conjugate to an
element in Kostant's section. Instead we consider  a maximal torus $T$ of~$G$
over~$k$ (not just over~$F$), so that in particular $T$ is unramified over~$F$,
and we assume that $T$ comes from a maximal torus $T_H$ in $H$.  We identify
$T$ and
$T_H$ using some admissible isomorphism, and for any $G$-regular
$\gamma \in \mathfrak t(F)$
we denote by $\gamma_H$ the corresponding element of $\mathfrak t_H(F)$.
It follows immediately from the definition of transfer factors on the
Lie algebra (see
\cite{waldspurger97} and perhaps also \cite{kottwitz98c}) that there
is a constant $c_T$
(depending only on the torus $T$ and our choice of admissible
isomorphism) such that for any
$G$-regular element $\gamma \in \mathfrak t(F)$
\[
\Delta(\gamma_H,\gamma)=c_T \cdot \Delta_{IV}(\gamma_H,\gamma)
\cdot\prod_\alpha(-1)^{\val(\alpha(\gamma))},
\]
where $\alpha$ runs over a set of representatives for the symmetric orbits of
Frobenius on the set of roots of $T$ in $G$ that do not come from
$H$. Here $\val$ is the
valuation on $\Fbar$ that takes the value $1$ on uniformizing
elements for~$F$. In fact
  the constant $c_T$ is equal to $1$, since any element
$\gamma \in \mathfrak t(\mathfrak o)$ whose image in $\mathfrak t(k)$
is $G$-regular (which implies
that $\val(\alpha(\gamma))=0$ for all roots $\alpha$) is
$G(F)$-conjugate to an element in our
choice of Kostant's section, and therefore for such $\gamma$ we have
(by \cite{kottwitz98c}) that
$\Delta(\gamma_H,\gamma)=\Delta_{IV}(\gamma_H,\gamma)$.

Now recall the homological transfer factor
$\Delta=\prod_{\alpha \in \Phi^+(G)-\Phi^+(H)}\partial_\alpha^{\val
\alpha(\gamma)}$. Frobenius acts on $X^*(T)$ by some element $\tau
\in \widetilde {\mathcal
W}^{G,H}_\gamma$, and we have $\tau(\Delta)=\eta(\tau) \cdot \Delta$
(see \S\ref{sec-endoscopic} for
the definitions of the group $\widetilde {\mathcal
W}^{G,H}_\gamma$ and the sign character $\eta$ on it). It is a simple
exercise to check that
\[
\eta(\tau)=\prod_\alpha(-1)^{\val(\alpha(\gamma))},
\]
with $\alpha$ again running over a set of representatives for the
symmetric orbits of
Frobenius on the set of roots of $T$ in $G$ that do not come from
$H$. Therefore we have the
equality
\begin{equation*}
\Delta(\gamma_H,\gamma)=\eta(\tau) \cdot \Delta_{IV}(\gamma_H,\gamma),
\end{equation*}
which for integral $\gamma$ can be rewritten as
\begin{equation}\label{transfer.factor}
\Delta(\gamma_H,\gamma)=\eta(\tau) \cdot q^{-r}
\end{equation}
with $r$ again defined  as $\dim(X_\gamma)-\dim(X^H_{\gamma_H})$.

\subsection{Remarks on the fundamental lemma for unramified
tori}\label{subsec-remarksFL}
We continue the discussion of the previous subsection.
Assume that the \'etale cohomology groups
$H^*(X_{\gamma}, \overline{\mathbb Q}_{\ell})$ and $H^*(X^H_{\gamma_H};
\overline{\mathbb Q}_{\ell})$ are pure.
By
\cite{GKM2} this is true in the equal valuation case, in other words,
when $\val(\alpha(\gamma))$ is
independent of
$\alpha$, since in that case there is a paving by affines defined over $k$.

We expect Theorems \ref{thm-cohomology} and \ref{thm-mainA} to have
corresponding statements in \'etale cohomology.  Assume this is so.
Then Theorem
\ref{thm-cohomology} will give a formula for the $T$-equivariant \'etale
homology of $X_{\gamma}$ as a quotient
\begin{equation*}\begin{CD}  0 @>>>
\sum_{\alpha \in \Phi^+_G} L_{\alpha, \gamma} @>>>
  \overline{\mathbb Q}_{\ell}[\Lambda] \otimes \bS(\mathfrak
t(\overline{\mathbb Q}_{\ell})) @>>> H_*^T(X_{\gamma},
\overline{\mathbb Q}_{\ell}) @>>> 0.  \end{CD}\end{equation*}
The surjection $\overline{\mathbb Q}_{\ell}[\Lambda] \otimes \bS(\mathfrak
t(\overline{\mathbb Q}_{\ell})) \to H_*^T(X_{\gamma},
\overline{\mathbb Q}_{\ell})$ is Frobenius equivariant, and the
action of Frobenius on
$\overline{\mathbb Q}_{\ell}[\Lambda] \otimes \bS(\mathfrak
t(\overline{\mathbb Q}_{\ell}))$ is given by $\tau \otimes Q$, where  $\tau
\in
\widetilde{\mathcal W}^{G,H}_{\gamma}$ is the element through which
the Frobenius element $\sigma$
acts on
$X_*(T)$ and
$Q$ denotes the endomorphism of the symmetric algebra induced by
multiplication by $q$ on the vector
space $\mathfrak t(\overline{\mathbb Q}_{\ell})$. It follows that
the action of the Frobenius element $\sigma$  on
$H_{2i}(X_\gamma;\overline{\mathbb Q}_{\ell})$
is given by
$\sigma =q^i \cdot \tau$.

  Therefore Corollary \ref{cor-ordinaryH} gives an
  isomorphism
\begin{equation*}
  H_*(X_{\gamma}; \overline{\mathbb Q}_{\ell})_J \simeq
H_*(X^H_{\gamma_H}; \overline{\mathbb
Q}_{\ell})_J[-2r](r)
\end{equation*}
that transforms by the sign $\eta(\tau)$ under the action of
Frobenius. As in the proof of Proposition \ref{prop-E2}, this in turn gives an
isomorphism
\begin{equation*}
\text{Tor}^{\Lambda}_p((\overline{\mathbb Q}_{\ell})_\kappa ,
H_q(X_{\gamma}; \overline{\mathbb
Q}_{\ell})
\simeq  \text{Tor}^{\Lambda}_p((\overline{\mathbb Q}_{\ell})_\kappa ,
H_{q-2r}(X^H_{\gamma_H}; \overline{\mathbb Q}_{\ell})(r).
\end{equation*}
Taking the vector space dual of both sides of this equation, we obtain an isomorphism
\begin{equation*}
\Ext^p_\Lambda((\overline{\mathbb Q}_{\ell})_\kappa, H^q(X_{\gamma};
\overline{\mathbb Q}_{\ell}))
\simeq  \Ext^p_\Lambda((\overline{\mathbb Q}_{\ell})_\kappa ,
H^{q-2r}(X^H_{\gamma_H}; \overline{\mathbb Q}_{\ell}))(-r),
\end{equation*}
that transforms by the sign $\eta(\tau)$ under the action of
Frobenius. It follows that
\begin{equation*}
q^{-r}\eta(\tau)\text{Tr}(\s^{-1};RHom^\bullet_{\Lambda}((\Qbar_l)
_\kappa,R\G(X_{\gamma},\Qbar_l)))=
\text{Tr}(\s^{-1};RHom^\bullet_{\Lambda}((\Qbar_l)_\kappa,R\G(X^H_{\gamma_H},\Qbar_l))).
\end{equation*}
Thus, under the assumptions
we have made, the fundamental lemma for unramified tori follows from
 this last equation
  together with  formulas \eqref{variant.main} and
\eqref{transfer.factor}.


\begin{thebibliography}{AMRT}

\bibitem[B96]{Bez} R. Bezrukavnikov, {\it The dimension of the fixed
point set on affine flag
manifolds}, Math. Res. Lett. {\bf 3} (1996), 185--189.

\bibitem[CS74]{CS} T. Chang and T. Skjelbred, {\it The topological Schur lemma
and related results}, Ann. Math. {\bf 100} (1974), 307--321.


\bibitem[D71]{HodgeII} P. Deligne, {\it Th\'eorie de Hodge II}, Publ. Math.
I.~H.~E.~S. {\bf 40} (1971), 5--57.


\bibitem[D75]{Deligne} P. Deligne, {\it Poids dans la cohomologie des
vari\'et\'es alg\'ebriques}, Proc. Int. Congr. Math. 1974 (Vancouver),
Canadian Math. Congress, 1975, pp. 79--85.

\bibitem[GKM98]{GKM1} M. Goresky, R. Kottwitz, and R. MacPherson, {\it Koszul
duality, equivariant cohomology, and the localization theorem}.  Invent. Math.
{\bf 131} (1998), 25--83.

\bibitem[GKMx]{GKM2} M. Goresky, R. Kottwitz, and R. MacPherson, {\it Purity
of equivalued affine Springer fibers}, to appear.

\bibitem[H86]{Hall}  M. Hall, {\bf Combinatorial Theory} (second edition),
John Wiley \& Sons, New York (1986).

\bibitem[Hum95]{Humphreys} J. Humphreys, {\bf Conjugacy Classes in Semisimple
Algebraic Groups}, Mathematical Surveys and Monographs {\bf 43}, Amer. Math.
Soc., Providence R.I., 1995.

\bibitem[KL80]{KL0} D. Kazhdan and G. Lusztig, {\it A topological approach to
Springer's representations}, Adv. in Math. {\bf 38} (1980), 222--228.

\bibitem[KL88]{KL}  D. Kazhdan and G. Lusztig, {\it Fixed point varieties on
affine flag manifolds}. Israel J. Math {\bf 62} (1988), 129--168.

\bibitem[Kot97]{Kottwitz} R. Kottwitz, {\it Isocrystals with additional
structure, II}, Comp. Math. {\bf 109} (1997), 255--339.

\bibitem[Kot99]{kottwitz98c} R. Kottwitz, {\it Transfer factors for
{L}ie algebras}, Represent.
Theory {\bf 3} (1999), 127--138.

\bibitem[Ku96]{Kumar} S. Kumar, {\it An introduction to {\it ind}-varieties},
in {\bf Vector Bundles On Curves--New Directions}, Lecture Notes in
Mathematics {\bf 1649}, Springer Verlag, N.Y., 1996, pp. 33--38.

\bibitem[La83]{Langlands} R. Langlands, {\bf Les d\'ebuts d'une formule des
traces stables}, monograph, Universit\'e de Paris VII, 1983.

\bibitem[LS87]{ls}  R. Langlands and D. Shelstad, {\it On the
definition of transfer factors},
Math. Ann. {\bf 278} (1987), 219--271.


\bibitem[L02a]{LaumonFL} G. Laumon, Sur le lemme fondamental pour les groupes 
unitaires, 
{\tt http:\slash\slash arxiv.org\slash math.AG\slash 0212245}

\bibitem[L02b]{Jacob} G. Laumon,  Fibres de Springer et Jacobiennes 
compactifiees, 
{\tt http:\slash\slash arxiv.org\slash math.AG\slash 0204109}



\bibitem[Lu96]{lusztig} G. Lusztig, {\it Affine Weyl groups and
conjugacy classes in Weyl groups}, Transform. Groups
{\bf 1} (1996), 83--97.



\bibitem[LS91]{LS}  G. Lusztig and J. M. Smelt, {\it Fixed point varieties on
the space of lattices}, Bull. London Math. Soc. {\bf 23} (1991), 213--218.



\quash{
\bibitem[Ro90]{Rogawski} J. Rogawski, {\bf Automorphic Representations of
Unitary Groups in Three Variables}, Annals of Math Studies {\bf 123},
Princeton University Press, Princeton NJ 1990. } 


\bibitem[Sa97]{Sage2} D. Sage, {\it A construction of affine Weyl group
representations} Comp. Math., {\bf 108} (1997), 241--245.

\bibitem[Sa00]{Sage} D. Sage, {\it The geometry of fixed point varieties on
affine flag manifolds}, Trans. Amer. Math. Soc. {\bf 352} (2000), 2087--2119.

\bibitem[Ser68]{serre68} J-P. Serre, {\bf Corps Locaux},  Hermann, Paris, 1968.

\bibitem[Sh82]{Shafarevich} I. Shafarevich, {\it On some infinite dimensional
groups II}, Math. USSR Izvestija {\bf 18} (1982), 185--194.

\bibitem[W97]{waldspurger97} J.-L. Waldspurger, {\it Le lemme
fondamental implique le transfert},
Compositio Math. {\bf 105} (1997), 153--236.

\end{thebibliography}
\end{document}